\def\LaTeX{%
  \let\Begin\begin
  \let\End\end
  \let\salta\relax
  \let\finqui\relax
  \let\futuro\relax}

\def\UK{\def\our{our}\let\sz s}
\def\USA{\def\our{or}\let\sz z}

\UK

\LaTeX

\USA

\salta

\documentclass[twoside,12pt]{article}
\usepackage[usenames,dvipsnames]{color}
\usepackage{amsmath}
\usepackage{amsthm}
\usepackage{amssymb}
\usepackage[mathcal]{euscript}

\usepackage[T1]{fontenc}
\usepackage[latin1]{inputenc}
\usepackage[english]{babel}
\usepackage[babel]{csquotes}

\usepackage{cite}

\usepackage{latexsym}
\usepackage{graphicx}
\usepackage{mathrsfs}
\usepackage{mathrsfs}
\usepackage{hyperref}
\usepackage{pgfplots}

\usepackage{amsfonts}
\usepackage{colortbl}
\usepackage{color}

\definecolor{viola}{rgb}{0.3,0,0.7}
\definecolor{ciclamino}{rgb}{0.5,0,0.5}

\def\pier #1{{\color{red}#1}}
\def\gianni #1{{\color{green}#1}}
\def\gabri #1{{\color{red}#1}}
\def\gabri  #1{#1}
\def\pier  #1{#1}
\def\gianni  #1{#1}

\def\oldpier #1{#1}

\def\oldrevis #1{#1}

\renewcommand{\theequation}{\thesection.\arabic{equation}}

\finqui

\def\Beq{\Begin{equation}}
\def\Eeq{\End{equation}}
\def\Bsist{\Begin{eqnarray}}
\def\Esist{\End{eqnarray}}

\def\Bthm{\Begin{theorem}}
\def\Ethm{\End{theorem}}

\def\Brem{\Begin{remark}\rm}
\def\Erem{\End{remark}}

\def\Bnot{\Begin{notation}\rm}
\def\Enot{\End{notation}}
\def\Bdim{\Begin{proof}}
\def\Edim{\QED\End{proof}}

\let\non\nonumber

\def\step #1 \par{\medskip\noindent{\bf #1.}\quad}

\def\Lip{Lip\-schitz}

\def\aand{\quad\hbox{and}\quad}

\def\lhs{left-hand side}
\def\rhs{right-hand side}

\def\multibold #1{\def\arg{#1}%
  \ifx\arg\pto \let\next\relax
  \else
  \def\next{\expandafter
    \def\csname #1#1#1\endcsname{{\bf #1}}%
    \multibold}%
  \fi \next}

\def\pto{.}

\def\multical #1{\def\arg{#1}%
  \ifx\arg\pto \let\next\relax
  \else
  \def\next{\expandafter
    \def\csname cal#1\endcsname{{\cal #1}}%
    \multical}%
  \fi \next}

\def\multimathop #1 {\def\arg{#1}%
  \ifx\arg\pto \let\next\relax
  \else
  \def\next{\expandafter
    \def\csname #1\endcsname{\mathop{\rm #1}\nolimits}%
    \multimathop}%
  \fi \next}

\multibold
qwertyuiopasdfghjklzxcvbnmQWERTYUIOPASDFGHJKLZXCVBNM.

\multical
QWERTYUIOPASDFGHJKLZXCVBNM.

\multimathop
diag dist div dom mean meas sign Sign supp .

\def\Span{\mathop{\rm span}\nolimits}

\def\accorpa #1#2{\eqref{#1}\pier{--}\eqref{#2}}
\def\Accorpa #1#2 #3 {\gdef #1{\eqref{#2}-\eqref{#3}}%
  \wlog{}\wlog{\string #1 -> #2 - #3}\wlog{}}

\def\separa{\noalign{\allowbreak}}

\def\somma #1#2#3{\sum_{#1=#2}^{#3}}

\def\graffe #1{\mathopen\{#1\mathclose\}}

\def\<#1>{\mathopen\langle #1\mathclose\rangle}
\def\norma #1{\mathopen \| #1\mathclose \|}

\def\[#1]{\mathopen\langle\!\langle #1\mathclose\rangle\!\rangle}

\def\iot {\int_0^t}

\def\intQt{\int_{Q_t}}
\def\intQ{\int_Q}
\def\iO{\int_\Omega}

\def\dt{\partial_t}
\def\dn{\partial_\nu}
\def\ds{\,ds}

\def\cpto{\,\cdot\,}

\def\checkmmode #1{\relax\ifmmode\hbox{#1}\else{#1}\fi}
\def\aeO{\checkmmode{a.e.\ in~$\Omega$}}
\def\aeQ{\checkmmode{a.e.\ in~$Q$}}

\def\aet{\checkmmode{a.e.\ in~$(0,T)$}}

\def\aat{\checkmmode{for a.a.~$t\in(0,T)$}}

\def\erre{{\mathbb{R}}}
\def\erren{\erre^n}

\def\genspazio #1#2#3#4#5{#1^{#2}(#5,#4;#3)}
\def\spazio #1#2#3{\genspazio {#1}{#2}{#3}T0}

\def\L {\spazio L}
\def\H {\spazio H}
\def\W {\spazio W}

\def\C #1#2{C^{#1}([0,T];#2)}

\def\Lx #1{L^{#1}(\Omega)}
\def\Hx #1{H^{#1}(\Omega)}

\def\HxG #1{H^{#1}(\Gamma)}

\def\LQ #1{L^{#1}(Q)}

\def\Luno{\Lx 1}
\def\Ldue{\Lx 2}
\def\Linfty{\Lx\infty}

\def\Huno{\Hx 1}
\def\Hdue{\Hx 2}
\def\Hunoz{{H^1_0(\Omega)}}

\def\LQ #1{L^{#1}(Q)}

\let\theta\vartheta
\let\eps\varepsilon
\let\badphi\phi
\let\phi\varphi

\let\TeXchi\chi                         
\newbox\chibox
\setbox0 \hbox{\mathsurround0pt $\TeXchi$}
\setbox\chibox \hbox{\raise\dp0 \box 0 }
\def\chi{\copy\chibox}

\long\def\salta #1\finqui{}

\def\QED{\hfill $\square$}

\let\hat\widehat

\def\normaV #1{\norma{#1}_V}
\def\normaH #1{\norma{#1}_H}

\def\normaVp #1{\norma{#1}_*}
\def\(#1\){\bigl(\!\bigl(#1\bigr)\!\bigr)}
\def\(#1\){\left(\!\left(#1\right)\!\right)}

\def\Beta{\widehat\beta}
\def\Pi{\widehat\pi}
\def\hatS{\hat S}
\def\phiz{\phi_0}
\def\chiz{\chi_0}
\def\phistar{\phi^*}
\def\rhostar{\rho^*}
\def\deltastar{\delta^*}
\def\Tstar{T^*}

\def\Csh{C_{sh}}
\def\Cstr{C_{str}}
\def\cO{C_\Omega}
\def\hatC{\hat C}
\def\mO{|\Omega|}
\def\mz{m_0}
\def\Vp{V^*}
\def\Vz{V_0}
\def\Vzp{V^*_0}
\def\Wz{W_0}
\def\Vmz{V_{mean}}
\def\Vmzp{V_{mean}^*}
\def\betaz{\beta^\circ}
\def\Sz{S^\circ}
\def\betastar{\beta^*}
\def\muG{\mu_\Gamma}
\def\muH{\mu_\calH}

\def\betaeps{\beta_\eps}
\def\Betaeps{\Beta_\eps}
\def\hatSeps{\hatS_\eps}
\def\hatseps{\hat\sigma_\eps}
\def\Seps{S_\eps}
\def\seps{\sigma_\eps}
\def\Feps{\calF_\eps}
\def\signeps{\sign\nolimits_\eps}
\def\phieps{\phi_\eps}
\def\mueps{\mu_\eps}

\def\chieps{\chi_\eps}
\def\Geps{G_\eps}
\def\weps{w_\eps}
\def\wz{w_0}
\def\Vm{V_m}
\def\Vn{V_n}
\def\Vinfty{V_\infty}
\def\phin{\phi^n}
\def\mun{\mu^n}
\def\phij{\badphi_j}
\def\etaj{\eta_j}
\def\beps{b^\eps}
\def\sigmaeps{\sigma^\eps}
\def\gn{g^n}

\def\soluz{(\phi,\mu,\xi,\zeta)}
\def\soluzeps{(\phieps,\mueps)}

\begin{document}

\pagestyle{myheadings}
\newcommand\testopari{\oldpier{\sc Colli \ --- \ Gilardi \ --- \ Marinoschi}}
\newcommand\testodispari{\oldpier{\sc \pier{Solvability and} sliding mode control \pier{for Cahn--Hilliard systems}}}
\markboth{\testopari}{\testodispari}

\begin{center}
{\Large \pier{Solvability and sliding mode control}\\[0.1cm]
for the viscous Cahn\pier{--}Hilliard system\\[0.2cm]
with a possibly singular potential}

\bigskip

\thispagestyle{empty}

Pierluigi Colli$^{1}$, Gianni Gilardi$^{1}$, Gabriela Marinoschi$^{2}$

\renewcommand*{\thefootnote}{\fnsymbol{footnote}}

\footnotetext{\textit{E-mail addresses}: 
\par
pierluigi.colli@unipv.it (P. Colli), \ 
\par
gianni.gilardi@unipv.it (G. Gilardi), 
\par
gabriela.marinoschi@acad.ro (G. Marinoschi)
\par
{}}\bigskip 

$^{1}$\oldpier{Dipartimento di Matematica ``F. Casorati'', Universit\`a di Pavia\\
and Research Associate at the IMATI -- C.N.R. Pavia\\ 
via Ferrata 5, 27100 Pavia, Italy}

$^{2}$\textquotedblleft Gheorghe Mihoc-Caius Iacob\textquotedblright\
Institute of Mathematical Statistics\\ and Applied Mathematics of the Romanian
Academy\\ 
Calea 13 Septembrie 13, 050711 Bucharest, Romania

\medskip
\end{center}

\begin{abstract}
In the present contribution we study a viscous Cahn--Hilliard system where a further
leading term in the expression for the chemical potential $ \mu$ is present. 
This term 
consists of a subdifferential operator \gianni{$S$ in~$L^2(\Omega)$ (where $\Omega$ is the domain where the evolution takes place)} 
acting on the difference of the
phase variable $\phi$ and a given state $\phistar $, which is prescribed and may depend on space and time. We prove
existence and continuous dependence results in case of \gianni{both homogeneous Neumann} and Dirichlet 
boundary conditions for the chemical potential~$\mu$. 
\pier{Next,} \gianni{by assuming that $S=\rho\sign$, a~multiple of the $\sign$ operator, and for smoother data, 
we first show regularity results. Then, in the case \pier{of} Dirichlet boundary conditions for~$\mu$ and 
under suitable conditions on $\rho$ and~$\Omega$, we also prove the sliding mode property, 
that is, that $\phi$~is forced 
to join the evolution of $\phistar $ in some time $T^*$  lower than the given final time~$T$}. 
We point out that all our results hold true for a very general and possibly singular multi-well potential acting on~$\varphi$.  
\vskip3mm

 \vskip3mm

\end{abstract}

\renewcommand{\theequation}{\arabic{section}.\arabic{equation}}


\section{Introduction}
\label{INTRO}
\setcounter{equation}{0}

This paper deals with the viscous Cahn--Hilliard system, which is further generalized
in order to admit an additional nonlinearity that plays as forcing term in order to reach a given evolution in the system. 
The resulting combination proposes an extension of the celebrated Cahn\pier{--}Hilliard system, 
a phenomenological model that has its origin in the work of J.W.~Cahn~\cite{Cahn1961}. 
In fact, Cahn studied the effects of interfacial energy on the stability of spinodal states in solid binary solutions, and \gianni{this}
took origin from a previous collaboration with J.W.~Hilliard \cite{Cahn1958}, where the functional
\begin{equation}
\label{eq:free_energy}
\mathcal F( \varphi)=\int_\Omega \left(\frac1 2 |\nabla \varphi |^2 + f(\varphi) \right)
\end{equation}
was proposed as a model for the free energy of a non-uniform system 
whose composition is described by the scalar field~$\varphi$. 
Actually, the viscous Cahn--Hilliard system considered here, which, up to our knowledge, 
was formulated by A.~Novick-Cohen~\cite{Novic1988viscous,NovicP1991TAMS}, reads
\Bsist
  && \dt\phi - \Delta\mu = 0
  \label{IprimaCH}
  \\
  && \tau \dt\phi - \Delta\phi + f'(\phi)
  = \mu + g
  \label{IsecondaCH}
\Esist
where the equations are understood to hold in a bounded domain $\Omega\subset\erre^3$
and in some time interval~$(0,T)$. 
Note that \eqref{IsecondaCH} contains the terms
$- \Delta\phi + f'(\phi)$ that can be interpreted as the variational derivative of~$\mathcal F( \varphi)$. 
The viscous contribution $\tau \dt\phi$ is completing the \lhs\ of~\eqref{IsecondaCH}.

In the above system, the variables $\phi$~and $\mu$ denote the order parameter and the associated chemical potential, respectively, with 
$\tau$ being the positive viscosity coefficient and $g$ standing for some \gabri{given} term.
In general, $f$~represents a non-convex potential; typical and physically significant examples for $\,f\,$ 
are the so-called {\em classical regular potential}, the {\em logarithmic double-well potential\/},
and the {\em double obstacle potential\/}, which are given, in this order,~by
\begin{align}
  & f_{reg}(r) := \frac 14 \, (r^2-1)^2 \,,
  \quad r \in \erre
  \label{regpot}
  \\
  & f_{log}(r) := \bigl( (1+r)\ln (1+r)+(1-r)\ln (1-r) \bigr) - c_1 r^2 \,,
  \quad r \in (-1,1)
  \label{logpot}
  \\[1mm]
  & f_{2obs}(r) :=  c_2 (1- r^2) 
  \quad \hbox{if $|r|\leq1$}
  \aand
  f_{2obs}(r) := +\infty
  \quad \hbox{if $|r|>1$}.
  \label{obspot}
\end{align}
Here, the constants $c_i$ in \eqref{logpot} and \eqref{obspot} satisfy
$c_1>1$ and $c_2>0$, so that $f_{log}$ and $f_{2obs}$ are nonconvex.
In cases like \eqref{obspot}, one has to split $f$ into a nondifferentiable convex part~$\Beta$ 
(the~indicator function of $[-1,1]$, in the present example) and a smooth perturbation~$\Pi$.
Accordingly, one has to replace the derivative of the convex part $\Beta$
by the subdifferential $\partial\Beta$ and interpret \eqref{IsecondaCH} as a differential inclusion. We will be more precise on that in the sequel. 

On the other hand, we have to emphasize that a wide number of generalizations of the
Cahn--Hilliard system have been proposed in the literature and it turns out that  these contributions are so many that it would be difficult to make a list here. 
We prefer to refer to the recent review paper\cite{Miranville2017}. 
In this respect, it is also worth 
pointing out that a systematic approach to derive and generalize the C-H system has been proposed by M.E.~Gurtin~\cite{Gurtin1996}, 
by extending the thermodynamical framework of Continuum Mechanics, as reported in \cite{Miranville2000} as well. 
We aim to mention an alternative procedure due to P.~Podio-Guidugli~\cite{Podio} that leads to another
viscous Cahn--Hilliard system of nonstandard type~\cite{CGPS1, CGPS2}. 
In addition, we report that in recent years Cahn--Hilliard and viscous Cahn--Hilliard systems have been 
employed successfully in many other branches of Science and Engineering, fields in which 
the segregation of a diffusant leads to pattern formation, such as population 
dynamics\cite{Liu2016}, image processing\cite{Bertozzi2007}, dynamics for mixtures 
of fluids\cite{GGW} and tumor modelling\cite{CGMR3}. 
In the case of the variable $\varphi$ understood as concentration, the recent paper~\cite{BCST} faces 
with a doubly nonlinear Cahn-Hilliard system, 
where both an internal constraint on the time derivative of $\varphi $ and the potential $f$ for $\varphi$ are introduced, 
thus leading to an equation more general than~\eqref{IsecondaCH}.
  
Coming back to the system~\accorpa{IprimaCH}{IsecondaCH},  we observe that proper supplementary conditions should complement it. 
As for~$\phi$, we consider the homogeneous Neumann boundary condition and the initial condition 
\Beq
  \dn\phi = 0
  \quad \hbox{on $\Gamma\times(0,T)$}
  \aand
  \phi(0) = \phiz
  \quad \hbox{in $\Omega$}
  \label{IbicCH}
\Eeq
where $\dn$ denotes the outward normal derivative on $\Gamma:=\partial\Omega$
and $\phiz$ is a given initial datum.
Regarding~$\mu$, in view of \eqref{IprimaCH} we have to add some boundary conditions.
We impose
\Beq
  \hbox{either} \quad
  \dn\mu = 0
  \quad \hbox{or} \quad
  \mu = \muG
  \quad \hbox{on $\Gamma\times(0,T)$}
  \label{Ibcmu}
\Eeq
where $\muG$ is a given boundary datum. 
We remark that the first condition of no flux through the boundary for $\mu$ is quite natural in the framework of Cahn--Hilliard systems 
and entails the mean value conservation property for~$\phi$, as the reader can easily realize 
by integrating \eqref{IprimaCH} over $\Omega \times (0,t)$, $t\in [0,T]$. 
Instead, the 
Dirichlet boundary condition for~$\mu$ (as in~\cite{WLFC08} and \cite{BCST}) is rather 
different and does not ensure any conservation, 
but it also looks reasonable from the modelling point of view and, as far as we know, 
classes of Dirichlet boundary data for $\mu$ are consistent with applications. 

In this paper, we study a sliding mode control (SMC) problem.
This consists in modifying the dynamics governed by the Cahn--Hilliard system 
by adding in equation \eqref{IsecondaCH} a further term
that forces the solution of the new system to satisfy $\phi(t)=\phistar(t)$
after some time~$\Tstar$, where $\phistar$ is a given function.
We stress that $\phistar$ is allowed to be time dependent, 
in contrast with the most part of the literature regarding sliding mode problems.
In fact, SMC is considered a classic instrument for regulation
of continuous or discrete systems in finite-dimensional settings
(see e.g.~the monographs \cite{BFPU08, EFF06, ES99, FMI11, I76, Utkin92, UGS09, YO99}),
in order to reach some stable states. 
\gabri{Here, we also want to allow} the possibility 
that the variable $\phi$ joins a prescribed evolution after the time $\Tstar$.

Then, the modified second equation is the following
\Beq
  \tau \dt\phi - \Delta\phi + f'(\phi) + S(\phi-\phistar)
  \ni \mu + g
  \label{Isecondagen}
\Eeq
in place of \eqref{IsecondaCH},
where now $S$ denotes some suitable maximal monotone graph in $\Ldue\times\Ldue$. 
In sliding mode problems, one generally has to choose operators that are singular at the origin.
So, a~typical choice for $S$ is given by the following~rule
\Beq
  \hbox{for $u,v\in\Ldue$}, \quad 
  v \in S(u)  
  \quad \hbox{means that} \quad
  v(x) \in \rho \sign u(x)
  \quad \aeO
  \label{IdefS}
\Eeq 
where $\rho$ is a positive parameter and $\sign$ is the subdifferential of the modulus $|\cdot|$, i.e.,
$\sign r:=r/|r|$ if $r\not=0$ and $\sign0:=[-1,1]$.
However, it is clear that the assumption that $S$ is a graph in $\Ldue\times\Ldue$
(not necessarily induced in $\Ldue$ by a graph in~$\erre\times\erre$)
is~much more general,
and the first result of ours regards well-posedness for such a generalized problem.
In proving it, we just have to reinforce the maximal monotonicity and subdifferential property 
by also assuming that $S$ grows at most linearly at infinity.
Moreover, we can manage both the Neumann and the Dirichlet boundary conditions given in~\eqref{Ibcmu}.

Next, under proper assumptions on the data of the problem that ensure 
further regularity for the solution,
we study the existence of a sliding mode.
This is done only for the Dirichlet boundary conditions for~$\mu$
and in the particular case~\eqref{IdefS}.
More precisely, $\rho>0$~has to be taken large enough
and we also have to assume the domain $\Omega$ to be small enough.
The result we prove is reminiscent of the one given in~\cite{CGMR3}, but here
we approach directly the viscous Cahn-Hilliard system and our key 
datum $\phistar$ is allowed to vary with time.

Now, we take the opportunity of reviewing some literature related to SMC, which offers a robust tool 
against abrupt variations, disturbances, time-delays, etc.\ in dynamics. 
The design procedure of a SMC scheme consists first in choosing a sliding set 
such that the original system restricted to it has a desired behavior, then 
\gianni{\pier{modifying} the dynamics in order to force} the involved variable to reach this set within a finite time. 
It is exactly for this aim that we add the term $ \rho \sign(\phi-\phistar)$ in the Cahn--Hilliard 
evolution for $\phi$ (cf.~\eqref{Isecondagen} and \eqref{IdefS}), in order to 
force $\phi$ to stay equal to a given desired \gianni{function} $\phistar$ in a finite time. 

Sliding mode controls are pretty interesting in applications and in recent years  
the extension of well-developed methods for finite-dimensional systems
(cf., e.g.,\cite{LO02, O83, O00, OU83}) and the control of infinite-dimensional 
dynamical systems (see~\cite{OU83, OU87, OU98}) have been faced. 
The theoretical development for PDE systems 
is still in its early stages: 
one can see the contributions~\cite{CRS11, GuWa,
Levaggi13, PPOU, XLGK13} dealing with semilinear PDE systems. 
We aim to quote~\cite{BCGMR}, where a sliding mode approach has been applied 
to phase field systems of Caginalp type: 
these systems combine  
the evolution of a phase variable to the one of the relative temperature, 
and the chosen SMC laws force the system to reach within finite time 
a sliding manifold. 
In that case it was possible to have different 
choices for the manifold:  in~\cite{BCGMR}, and also in~\cite{CM} 
which considers an extension of the \gianni{Caginalp} model, either one of the physical 
variables or a combination of them could reach a stable state. 
With reference to the results of \cite{BCGMR, CM}, we mention 
the analyses developed in \cite{Colt1, Colt2}: 
in particular, the second contribution is devoted to a conserved phase field system 
with a SMC feedback law for the internal energy in the temperature equation.

An outline of the present paper is as follows. 
In Section~\ref{STATEMENT} we state
precisely the problem, making clear the assumptions and presenting the different results we are going to prove. 
Section~\ref{UNIQUENESS} brings the proof of the continuous 
dependence result, which ensures uniqueness at least for the component variable~$\varphi$.
The approximation of the problem, based on Yosida regularizations of graphs and a Faedo--Galerkin scheme, 
is discussed in Section~\ref{APPROXIMATION}. 
The existence of solutions 
is shown in Section~\ref{EXISTENCE} by proving some a priori estimates and passing to the 
limit with respect to the parameter of the Yosida regularizations. 
Finally, Section~\ref{SLIDING} is completely devoted to the proof of the sliding mode property, 
first dealing with the regularity of the solution, then proving the existence of a 
suitable time~$\Tstar$, 
after which it occurs that \gianni{$\phi(t)=\phistar(t)$}. 


\section{Statement of the problem}
\label{STATEMENT}
\setcounter{equation}{0}

As in the Introduction, $\Omega$~is the domain where the evolution process takes place.
We assume that $\Omega$ 
is a bounded and connected open set in~$\erre^3$
(more generally, one could take $\Omega\subset\erre^d$ with $1\leq d\leq 3$),
which is supposed to have a smooth boundary~$\Gamma:=\partial\Omega$,
and we write $\mO$ and $\dn$ for the volume of~$\Omega$ 
and the outward normal derivative on~$\Gamma$, respectively.
Given some final time~$T>0$, we set for convenience
\Beq
  Q_t := \Omega \times (0,t)
  \quad \hbox{for $t\in(0,T]$}
  \aand
  Q := Q_T \,.
  \label{defQt}
\Eeq
If $X$ is a Banach space, $\norma\cpto_X$ denotes both its norm and the norm of~$X^3$. 
Moreover, the dual space of $X$ and the dual pairing between $X^*$ and~$X$
are denoted by $X^*$ and $\<\cpto,\cpto>$, the latter without indeces
since the choice of the space $X$ is clear every time from the context.
The only exception from the convention for the norms is given
by the the spaces $L^p$ constructed on $\Omega$ and~$Q$ for $p\in[1,\infty]$, 
whose norms are denoted by~$\norma\cpto_p$. 
Furthermore, we~put
\Bsist
  && H := \Ldue \,, \quad  
  V := \Huno 
  \aand
  W := \graffe{v\in\Hdue:\ \dn v=0}
  \label{defspazi}
  \\
  && \Vz := \Hunoz
  \aand
  \Wz := \Hdue \cap \Hunoz \,.
  \label{defspaziz}
\Esist
We endow these spaces with their standard norms.
Moreover, we identify $H$ with a subspace of $\Vp$ in the usual way, i.e.,
in order that $\<u,v>=\iO uv$ for every $u\in H$ and $v\in V$,
and obtain the Hilbert triplet $(V,H,\Vp)$.
Analogously, we consider the Hilbert triplet $(\Vz,H,\Vzp)$ when dealing with Dirichlet boundary conditions.

Now, we list our assumptions on the structure of the system at once.
We assume~that
\Bsist
  && \hbox{$\tau$ is positive real number}
  \label{hptau}
  \\
  && \beta := \partial\Beta, \quad
  S := \partial\hatS
  \aand
  \pi := \Pi{}'
  \quad \hbox{where}
  \label{defbetaApi}
  \\
  && \Beta : \erre \to [0,+\infty]
  \quad \hbox{is convex, proper and l.s.c.\ with}
  \quad \Beta(0) = 0
  \label{hpBeta}
  \\
  && \hatS : H \to \erre
  \quad \hbox{is convex, proper\pier{, l.s.c.}\ and $S$ satisfies} \enskip
  \normaH v \leq C_S (\normaH u + 1)
  \hskip 1.5em
  \non
  \\
  && \quad \hbox{for some constant $C_S$ and every $u\in H$ and $v\in S(u)$}
  \label{hpS}
  \\
  && \Pi : \erre \to \erre
  \quad \hbox{is of class $C^1$ with a \Lip\ continuous first derivative}.
  \qquad
  \label{hpPi}
\Esist
\Accorpa\HPstruttura hptau hpPi
Notice that all of the important examples \accorpa{regpot}{obspot} satisfy the above assumptions.
We also remark that $\beta$ and $S$ are maximal monotone graphs.
We denote by $D(\beta)$ the effective \pier{domain} of $\beta$ and, for $r\in D(\beta)$, 
we use the notation $\betaz(r)$ for the element of $\beta(r)$ having minimum modulus.
The similar notation $\Sz(u)$ for $u\in H$ refers to the minimum norm. 
For simplicity, we still write $\beta$ and $S$ for the graphs induced 
in $\Ldue$ and~$\LQ2$ by $\beta$ and in $\LQ2$ by~$S$, respectively.

As for the data of the problem we assume that
\Bsist
  && g \in \L2H
  \label{hpg}
  \\
  && \phiz \in V 
  \aand
  \Beta(\phiz) \in \Luno
  \label{hpphiz}
  \\
  && \phistar \in \L2H \,.
  \label{hpphistar}
\Esist
\Accorpa\HPdati hpg hpphistar
Moreover, in the case of the Neumann boundary conditions for~$\mu$, we also assume that
\Beq
  \hbox{$\mz := \mean\phiz$ belongs to the interior of $D(\beta)$}
  \label{hpmz}
\Eeq
where the symbol $\mean v$ denotes the mean value of the generic function $v\in\Luno$.
More generally (by~denoting by $1$ the function that is identically $1$ on~$\Omega$), we~set
\Beq
  \mean v : = \frac 1 {\mO} \, \< v , 1 >
  \quad \hbox{for every $v\in\Vp$}
  \label{defmean}
\Eeq
and it is clear that $\mean v$ is the usual mean value of $v$ if $v\in H$.

At this point, we can state the problem given by equations \eqref{IprimaCH} and \eqref{Isecondagen}
and the boundary and initial conditions given in \eqref{IbicCH} and~\eqref{Ibcmu}.
We first distinguish between the two different boundary conditions for~$\mu$.
Then, we unify the two problems.
We start with the case of the Neumann boundary conditions.

\step
The case of the Neumann boundary conditions

We write a variational formulation.
We set $\calV:=V$ for convenience and look for a quadruple $\soluz$ satisfying
\Bsist
  && \phi \in \H1H \cap \L2 V
  \label{regphi}
  \\
  && \mu \in \L2\calV
  \label{regmu}
  \\
  && \xi \in \L2H
  \aand 
  \xi \in \beta(\phi) \quad \aeQ
  \label{regxi}
  \\
  && \zeta \in \L2H
  \aand
  \zeta(t) \in S(\phi(t)-\phistar(t))
  \quad \aat
  \qquad
  \label{regzeta}
\Esist
\Accorpa\Regsoluz regphi regzeta
and solving the following system
\Bsist
  && \iO \dt\phi(t) \, v 
  + \iO \nabla\mu(t) \cdot \nabla v
  = 0 
  \quad \hbox{\aat\ and every $v\in\calV$}
  \qquad
  \label{prima}
  \\
  && \tau \iO \dt\phi(t) \, v
  + \iO \nabla\phi(t) \cdot \nabla v
  + \iO \bigl( \xi(t) + \pi(\phi(t)) + \zeta(t) \bigr) v
  \non
  \\
  && = \iO \bigl( \mu(t) + g(t) \bigr) v
  \quad \hbox{\aat\ and every $v\in V$}
  \label{seconda}
  \\
  && \phi(0) = \phiz \,.
  \label{cauchy}
\Esist
\Accorpa\Pbl prima cauchy

\Brem
\label{Piureg}
In fact, every solution enjoys some more regularity,
namely
\Beq
  \phi\in \L2W
  \aand
  \mu \in \L2W
  \label{piureg}
\Eeq
so that equations \eqref{prima} and \eqref{seconda} can be written
in the strong form \eqref{IprimaCH} and \eqref{Isecondagen}
complemented with homogeneous Neumann boundary conditions.
Indeed, both \eqref{prima} and \eqref{seconda} have the form
(with $u=\mu$ and $u=\phi$, respectively)
\Beq
  \iO \nabla u(t) \cdot \nabla v
  = \iO \psi(t) v
  \quad \hbox{\aat\ and every $v\in V$}
  \non
\Eeq
with $\psi\in\L2H$.
This implies that 
\Beq
  u \in \L2W
  \aand
  - \Delta u = \psi
  \quad \aeQ .
  \label{pde}
\Eeq
However, in connection with the general result stated below,
we just deal with the variational formulation
and do not need such a further regularity of the solution.
On the contrary, this remark is used in the last section.
\Erem

\Brem\label{Meanvalue}
It is worth noting that every solution also satisfies
\Beq
  \iO \dt\phi = 0 \quad \aet , 
  \quad \hbox{i.e.,} \quad
  \mean\phi(t) = \mean\phiz
  \quad \hbox{for every $t\in[0,T]$}
  \label{meanvalue}
\Eeq
as one immediately sees by choosing $v=1\in\calV=V$ in~\eqref{prima}.
\Erem

\step
The case of the Dirichlet boundary conditions

In the case of the Dirichlet boundary conditions for $\mu$ given in~\eqref{Ibcmu},
we first make the basic assumption on~$\muG$.
Since we still look for $\mu$ in $\L2\Huno$, we require~that
\Beq
  \muG \in \L2{\HxG{1/2}} \,.
  \label{hpmuG}
\Eeq
As for the problem, we have to force $\mu=\muG$, explicitly, 
and modify \eqref{prima} as far as the test functions are concerned.
Namely, \eqref{prima} is required to hold just for $v\in\Vz$.
However, it is convenient to reduce the boundary condition $\mu=\muG$ to the homogeneous one.
This can be done by introducing the harmonic extension $\muH$ of~$\muG$,
which is defined \aat~by
\Beq
  \Delta\muH(t) = 0 
  \quad \hbox{in $\Omega$} 
  \aand
  \muH(t) = \muG(t)
  \quad \hbox{on $\Gamma$} 
  \label{harmonic}
\Eeq
and by considering the problem that the difference $\mu-\muH$ has to solve.
However, it is better to avoid a new notation (as~we see in a moment)
and still term $\mu$ the above difference.
Then, for the new~$\mu$, both the regularity requirement \eqref{regmu} and the first equation \eqref{prima} 
remain unchanged provided that we set $\calV=\Vz$ now,
while the forcing term $g$ in \eqref{seconda} has to be replaced by the difference $g_*:=g-\muH$.
Hence, \eqref{seconda} formally remains unchanged too 
provided that we still use the symbol $g$ for the difference~$g_*$.
Notice that the new $g$ belongs to $\L2H$ as the old one (see \eqref{hpg})
since \eqref{hpmuG} trivially implies $\muH\in\L2H$.
Therefore the two problems corresponding to the two different boundary conditions for $\mu$ are unified
and we just have different meanings of $\mu$ and~$g$ in the two cases.

\Brem
\label{PiuregD}
More regularity for the old chemical potential $\mu$ as in Remark~\ref{Piureg} 
is ensured whenever $\muG\in\L2{\HxG{3/2}}$.
Indeed, this implies that $\muH\in\L2\Hdue$.
On the contrary, the mass conservation stated in Remark~\ref{Meanvalue}
cannot be expected in the case of the Dirichlet boundary conditions for~$\mu$.
Indeed, the choice $v=1$ in \eqref{prima} is no longer allowed since $\calV=\Vz$ now.
\Erem
  
\Bnot
\label{Notation}
From now on, it is understood that $g$ has its new meaning 
in the case of the Dirichlet boundary conditions for $\mu$
and the Dirichlet datum $\muG$ and its harmonic $\muH$ are not mentioned any longer in the problem.
The only exception to this rule will happen if we need more regularity for~$g$.
In that case, we come back to this point 
and give sufficient conditions on $\muG$ in order to satisfy the new requirements.

Furthermore, in the sequel of the paper, we simply write $\calV=V$ and $\calV=\Vz$
to identify the boundary conditions for~$\mu$
of the Neumann and Dirichlet type, respectively, in problem \Pbl.
\Enot

If $\calV=\Vz$, the component $\mu$ of any solution is uniquely determined 
whenever uniqueness holds for the first component~$\phi$,
since the first equation \eqref{prima} is uniquely solvable for $\mu$ in this case.
On the contrary, if $\calV=V$,
it is clear that no uniqueness for $\mu$ can be expected unless both $\beta$ and $S$ are single-valued,
and this is not the case in this paper.
Hence, the best one can have is just existence of a solution $\soluz$ 
and uniqueness and some continuous dependence for the first component of the solution
and uniqueness for the second component if $\calV=\Vz$.
In the next sections, the following result is proved:

\Bthm
\label{Wellposedness}
Let the assumptions \HPstruttura\ on the structure and \HPdati\ on the data be satisfied.
In addition, assume either $\calV=V$ and \eqref{hpmz} or $\calV=\Vz$.
Then, there exists a quadruplet $\soluz$ satisfying \Regsoluz\ and solving problem \Pbl.
Moreover, the component $\phi$ of any solution is uniquely determined in any case
and $\mu$ is uniquely determined if $\calV=\Vz$.
Furthermore, let $g_i$, $\phi_{0,i}$ and $\phistar_i$, $i=1,2$ be two choices of the data
and assume that $\phi_{0,1}$ and $\phi_{0,2}$ have the same mean value if $\calV=V$.
Then, for the first components of any corresponding solutions $(\phi_i,\mu_i,\xi_i,\zeta_i)$,
the continuous dependence inequality
\Bsist
  && \norma{\phi_1-\phi_2}_{\L\infty H\cap\L2V}
  \non
  \\
  && \leq C_{cd} \bigl(
    \norma{g_1-g_2}_{\L2H}
    + \normaH{\phi_{0,1}-\phi_{0,2}} 
    + \norma{\phistar_1-\phistar_2}_{\L2H}^{1/2}
  \bigr)
  \label{contdep}
\Esist
holds true for some constant $C_{cd}$ only depending on the structure of the system, $\Omega$, $T$
and an upper bound $M$ for the norms of $\phi_i$ and $\phistar_i$ in~$\L2H$.
\Ethm

The second result of the present paper, whose proof is given in Section~\ref{SLIDING}, is the existence of a sliding mode.
We cannot treat the Neumann boundary conditions for~$\mu$, unfortunately,
and just consider the case of the Dirichlet boundary conditions.
Moreover, the convex function $\hatS$ and its subdifferential $S$ have a particular shape, as said in the Introduction.
Precisely, we assume that
\Beq
  \hatS(u) = \rho \iO |u| 
  \quad \hbox{for $u\in H$}
  \label{hphatS}
\Eeq
where $\rho$ is a positive real number.
Then, $S$~is given by~\eqref{IdefS}, i.e., it is the graph induced in $H$ by $\rho\sign$, 
where we recall that $\sign$, the subdifferential of the modulus, is given~by
\Beq
  \sign r := \frac r {|r|}
  \quad \hbox{if $r\not=0$}
  \aand
  \sign 0 = [-1,1] .
  \label{defsign}
\Eeq
Moreover, our hypotheses on the data have to be reinforced.
Besides \HPdati, we require that
\Bsist
  && g \in \H1H \cap \LQ\infty
  \label{hpgbis}
  \\
  && \phiz \in W 
  \aand
  \betaz(\phiz) \in \Linfty
  \label{hpphizbis}
  \\
  && \phistar \in \W{2,1}\Luno \cap \L2W, \quad
  \phistar,\, \dt\phistar ,\, \Delta\phistar ,\, \betaz(\phistar) \in \LQ\infty \,.
  \qquad
  \label{hpphistarbis}
\Esist
\Accorpa\HPdatibis hpgbis hpphistarbis

\Brem
\label{PiuregmuG}
According to Notation~\ref{Notation},
the function $g$ appearing in \eqref{hpgbis} is the difference between the original forcing term
and the harmonic extension $\muH$ of the inhomogeneous boundary datum~$\muG$.
In order to satisfy~\eqref{hpgbis}, we have to assume the same regularity for both
the original forcing term and~$\muH$.
To obtain a sufficient condition for the latter, one can reinforce \eqref{hpmuG}
by also assuming that
\Beq
  \muG \in L^\infty(\Gamma\times(0,T))
  \aand
  \dt\muG \in L^2(\Gamma\times(0,T)).
  \label{hpmuGbis}
\Eeq
\Erem

\Brem
\label{FurhterHPdata}
Our new assumptions on $\phistar$ are just sufficient conditions for the existence of a sliding mode.
However, some of them are necessary since they are satisfied by~$\phi$.
Indeed, as the other data are more regular, more regularity for $\phi$ is expected.
Moreover, we observe that, in the case of an everywhere defined potential like~\eqref{regpot},
the boundedness condition on $\betaz(\phistar)$ is satisfied whenever $\phistar$ is bounded.
On the contrary, when dealing with potentials like \eqref{logpot} or~\eqref{obspot},
we also have to require the smallness condition $\norma\phistar_\infty<1$.
\Erem

As announced in the Introduction, we also have to assume that $\rho$ is large enough and $\Omega$ is small enough.
More precisely, once we fix the class $\calO$ of the domains of $\erre^3$ 
that have the same shape of~$\Omega$,
then $\mO$~has to be small enough.
To better explain this condition, we fix a class $\calO$ as said above.
Then, there exists a constant $\Csh$ realizing the inequalities
\Bsist
  && \norma v_\infty
  \leq \Csh \mO^{1/6} \normaH{\Delta v}
  \quad \hbox{for every $v\in\Wz$} 
  \label{embeddingz}
  \\
  && \norma v_\infty
  \leq \Csh \bigl( \mO^{-1/2} \normaH v + \mO^{1/6} \normaH{\Delta v} \bigr)
  \quad \hbox{for every $v\in W$} 
  \label{embedding}
\Esist
whenever $\Omega\in\calO$.
The smallness condition on $\mO$ will involve the constant~$\Csh$.

\Brem
\label{Cshape}
We summarize the argument of \cite[Rem.~2.1]{CGMR3}
that shows that the constant $\Csh$ realizing \accorpa{embeddingz}{embedding} actually exists.
The prototype for $\Omega$ in the given class $\calO$
is an open set $\Omega_0\subset\erre^3$
(which \pier{is} supposed to be bounded, connected and smooth)
with $|\Omega_0|=1$,
and the general set $\Omega\in\calO$ has the form
\Beq  
  \Omega = x_0 + \lambda R\,\Omega_0
  \label{hpOmega}
\Eeq
where $x_0$ is a point in $\erre^3$, 
the real number $\lambda$ is positive and $R$ belongs to the the rotation group $SO(3)$.
We first notice that our assumptions on~$\Omega_0$, the continuous embedding $W\subset\Linfty$ and elliptic regularity
ensure that the inequalities \eqref{embeddingz} and \eqref{embedding} hold true for some constant~$\Csh$
(depending only on~$\Omega_0$) 
if~$\Omega$ and $\mO$ are replaced replaced by $\Omega_0$ and~$1$, respectively.
Now, assume that the domain $\Omega$ we are dealing is given by \eqref{hpOmega} as said above.
Then, it is easy to check that $\mO=\lambda^3$, i.e., $\lambda=\mO^{1/3}$,
and that \accorpa{embeddingz}{embedding} are still satisfied for $\Omega$ with the same constant~$\Csh\,$.
\Erem

Here is our \gabri{result}, which provides both further regularity for the solution
and the existence of a sliding mode.

\Bthm
\label{Sliding}
In addition to the hypotheses of Theorem~\ref{Wellposedness},
assume $\calV=\Vz$, \eqref{hphatS} on the function $\hatS$ and \HPdatibis\ on the data.
Then, every solution $\soluz$ to problem \Pbl\ enjoys the further regularity
\Bsist
  && \phi \in \W{1,\infty}H \cap \H1V \cap \L\infty W
  \subset \LQ\infty
  \label{regphibis}
  \\
  && \mu \in \L\infty W
  \subset \LQ\infty
  \label{regmubis}
  \\
  && \xi \in \L\infty H \,.
  \label{regxibis}
\Esist
\Accorpa\Regsoluzbis regphibis regxibis
Moreover, assume that $\Omega$ belongs to a class $\calO$ of open sets 
for which \eqref{embedding} is guaranteed.
Then, there exist $\rhostar>0$ and $\deltastar>0$ such that the following holds true:
if $\rho>\rhostar$ and $\mO<\deltastar$, then the component $\phi$ of any solution
satisfies for some $\Tstar\in(0,T)$ the sliding condition $\phi(t)=\phistar(t)$ for every $t\in[\Tstar,T]$.
\Ethm

\Brem
\label{Remregsoluz}
The properties specified in the above statement refer to any solution.
As for the last sentence, we recall that the component $\phi$ of any solution $\soluz$ is uniquely determined
(as~well as $\mu$ since $\calV=\Vz$).
On the contrary, no uniqueness for the components $\xi$ and $\zeta$ is ensured since $S$ is multivalued.
Nevertheless, the regularity property \eqref{regxibis} holds for every solution.
Indeed, \eqref{seconda}~can be written as a PDE \aeQ\ (apply Remark~\ref{Piureg} to this equation),
so that \eqref{regxibis} follows by comparison, since all the other \pier{terms} of the PDE belong to~$\L\infty H$,
due~to \accorpa{regphibis}{regmubis} and the boundedness of the (possibly non unique) component~$\zeta$, 
since $S$ is the graph induced in $H$ by $\rho\sign$.
\Erem

The rest of the paper is organized as follows.
In the next section, we prove the part of Theorem~\ref{Wellposedness} concerning uniqueness and continuous dependence.
The existence part is concluded in Section~\ref{EXISTENCE} and prepared in Section~\ref{APPROXIMATION},
where an approximating problem is introduced and solved.
The proof of Theorem~\ref{Sliding} is presented in the last Section~\ref{SLIDING}.

In proving our results, we make a wide use of the Schwarz and Young inequalities.
We recall the latter:
\Beq
  ab \leq \delta a^2 + \frac 1{4\delta} \, b^2
  \quad \hbox{for every $a,b\in\erre$ and $\delta>0$}.
  \label{young}
\Eeq
Moreover, we often account for the Poincar\'e inequalities
\Bsist
  && \normaV v \leq \cO \, \bigl( \normaH{\nabla v} + |\mean v| \bigr)
  \aand
  \normaV v \leq \cO \, \normaH{\nabla v}
  \non
  \\
  && \hbox{for every $v\in V$ and every $v\in\Vz$, respectively,}
  \label{poincare}
\Esist
with a constant $\cO$ that only depends on~$\Omega$.
Furthermore, we take advantage of a tool that is rather common 
in the study of problems related to the Cahn--Hilliard equations.
In order to introduce it, we recall that $\Omega$ is connected,
define the subspaces 
\Beq
  \Vmz := \graffe{v\in V:\ \mean v=0}
  \aand
  \Vmzp := \graffe{v\in\Vp:\ \mean v=0}
  \label{defVmzVmzp}
\Eeq
and consider, for $\psi\in\Vp$, the problem of finding
\Beq
  u \in V
  \quad \hbox{such that} \quad
  \iO \nabla u \cdot \nabla v
  = \< \psi , v >
  \quad \hbox{for every $v\in V$}.
  \label{neumann}
\Eeq
By the way, if $\psi\in H$, this is the usual Neumann problem
\Beq
  - \Delta u = \psi
  \quad \hbox{in $\Omega$}
  \aand
  \dn u = 0
  \quad \hbox{on $\Gamma$}.
  \non
\Eeq
Now, for $\psi\in\Vp$, \eqref{neumann} is solvable if and only if $\psi\in\Vmzp$.
Moreover, if $\psi\in\Vmzp$, exactly one of the solutions belongs to~$\Vmz$.
This implies that the operator
\Bsist
  && \calN: \Vmzp \to \Vmz
  \quad \hbox{defined by the following rule:}
  \non
  \\
  && \hbox{for $\psi\in\Vmzp$,\quad $\calN\psi$ is the unique solution $u$ to \eqref{neumann} belonging to $\Vmz$}
  \qquad\qquad
  \label{defN}
\Esist
is well defined.
It turns out that $\calN$ is an isomorphism and that the function
\Bsist
  && \Vp \ni \psi \mapsto \normaVp\psi^2
  := \normaH{\nabla\calN(\psi-\mean\psi)}^2 + |\mean\psi|^2
  \non
  \\
  && \quad \bigl( \hbox{in particular,\quad $\normaVp\psi^2=\normaH{\nabla\calN\psi}^2=\<\psi,\calN\psi>$\quad if $\psi\in\Vmzp$} \bigr)
  \label{normaVp}
\Esist
is the square of a norm on $\Vp$ that is equivalent to the standard one.
In the sequel, we use \eqref{normaVp} for the norm in~$\Vp$.
Similarly, we introduce the Dirichlet problem solver related to the Poisson equation
\Beq
  - \Delta u = \psi
  \quad \hbox{or} \quad
  \iO \nabla u \cdot \nabla v
  = \< \psi , v >
  \quad \hbox{for every $v\in\Vz$}
  \label{poisson}
\Eeq
with homogeneous Dirichlet boundary conditions.
It is the operator
\Bsist
  && \calD: \Vzp \to \Vz
  \quad \hbox{defined by the following rule:}
  \non
  \\
  && \hbox{for $\psi\in\Vzp$,\quad $\calD\psi$ is the unique solution $u$ to \eqref{poisson} belonging to $\Vz$}.
  \qquad\qquad
  \label{defD}
\Esist
We notice that the function
\Beq
  \Vzp \ni \psi \mapsto \normaVp\psi
  := \normaH{\nabla\calD\psi}
  \label{normaVzp}
\Eeq
is a norm on $\Vzp$ that is equivalent to the standard one since $\calD$ is an isomorphism. 
We remark that
\begin{align}
  & \< \dt v(t) , \calL v(t) >
  = \< v(t) , \calL(\dt v(t)) >
  = \frac 12 \, \frac d{dt} \, \normaVp{v(t)}^2
  \quad \aat
  \non
  \\
  & \quad \hbox{for every $v\in\H1{\calV^*}$, where $\calL=\calN$ if $\calV=V$ and $\calL=\calD$ if $\calV=\Vz$} \,.
  \label{propND} 
\end{align}
In \eqref{propND}, the notation $\normaVp\cpto$ 
means \eqref{normaVp} if $\calV=V$ and \eqref{normaVzp} if $\calV=\Vz$, of course.
Also in the sequel, the meaning of $\normaVp\cpto$ is clear from the context and no confusion can arise.


\section{Partial uniqueness and continuous dependence}
\label{UNIQUENESS}
\setcounter{equation}{0}

In this section, we prove the part of Theorem~\ref{Wellposedness} regarding partial uniqueness and continuous dependence.
Namely, we just prove the latter, since the former follows as a consequence.
As for uniqueness of $\mu$ in the case $\calV=\Vz$,
we have already noticed that the first equation \eqref{prima} is uniquely solvable for $\mu$ in this case,
so that uniqueness of $\mu$ follows from uniqueness for~$\phi$.
So, we pick two choices $g_i$, $\phi_{0,i}$ and $\phistar_i$, $i=1,2$, of the data and a constant $M$ as in the statement.
We assume that $(\phi_i,\mu_i,\xi_i,\zeta_i)$ are arbitrary corresponding solutions
and we prove the continuous dependence inequality~\eqref{contdep}.
For brevity, we use the same symbol~$c$ 
(even in the same line or a chain of inequalities)
for different constants depending only on the structure of our system, $\Omega$, $T$ and~$M$.
We set for convenience
\Bsist
  && g := g_1 - g_2 \,, \quad
  \phiz := \phi_{0,1} - \phi_{0,2} \,, \quad
  \phistar := \phistar_1 - \phistar_2
  \non
  \\
  && \phi := \phi_1 - \phi_2 \,, \quad
  \mu := \mu_1 - \mu_2 \,, \quad
  \xi := \xi_1 - \xi_2
  \aand
  \zeta := \zeta_1 - \zeta_2 \,.
  \non
\Esist
We recall that $\mean\phi_{0,1}=\mean\phi_{0,2}$ if $\calV=V$. 
In this case, by applying \eqref{meanvalue} to both $\phi_1$ and $\phi_2$,
we infer that $\mean\phi$ vanishes identically.
Hence, we can write \eqref{prima} at the time $s$ for both \pier{solutions}
and take $v=\calN\phi(s)$ as test function in the difference.
In the case $\calV=\Vz$ of the Dirichlet boundary conditions,
we can take $v=\calD\phi(s)$ without any trouble.
In order to unify the arguments,
we set $\calL=\calN$ and $\calL=\calD$ according to $\calV=V$ or $\calV=\Vz$.
We recall that the symbol $\normaVp\cpto$ means \eqref{normaVp} or \eqref{normaVzp}
according to $\calV=V$ or $\calV=\Vz$.
Thus, we have in both cases
\Beq
  \iO \dt\phi(s) \, \calL\phi(s)
  + \iO \nabla\mu(s) \cdot \nabla\calL\phi(s)
  = 0 \,.
  \non
\Eeq
By then integrating over $(0,t)$ with $t\in(0,T)$ and recalling~\eqref{propND}, 
we~have
\Beq
  \frac 12 \, \normaVp{\phi(t)}^2
  + \intQt \nabla\mu \cdot \nabla\calL\phi
  = \frac 12 \, \normaVp\phiz^2 \,.
  \non
\Eeq
Similarly, we write \eqref{seconda} for both solutions at the time~$s$,
test the difference by $\phi(s)$ and integrate over~$(0,t)$.
Moreover, we subtract the same quantity $\intQt\zeta\phistar$ to both sides and rearrange.
We obtain
\Bsist
  && \frac \tau 2 \, \normaH{\phi(t)}^2
  + \intQt |\nabla\phi|^2
  + \intQt \xi \phi
  + \intQt \zeta (\phi-\phistar)
  \non
  \\
  && = \frac \tau 2 \, \normaH\phiz^2
  + \intQt \mu \phi
  + \intQt g \phi
  - \intQt \bigl( \pi(\phi_1) - \pi(\phi_2) \bigr) \, \phi
  - \intQt \zeta \phistar .
  \non
\Esist
At this point, we add the above equalities to each other.
Then, the terms involving $\mu$ cancel out by the definition of~$\calL$
(see \eqref{defN} if $\calL=\calN$ and \eqref{defD} if $\calL=\calD$)
and those containing $\xi$ and $\zeta$ on the \lhs\ are nonnegative by monotonicity.
Furthermore, we can owe to the \Lip\ continuity of $\pi$ given by \eqref{hpPi} 
and estimate the corresponding term.
As for the last integral on the \rhs, we can account for the linear growth of $S$ given by \eqref{hpS} 
(since $\zeta_i(s)\in S(\phi_i(s)-\phistar_i(s))$ for $i=1,2$)
and use the upper bound $M$ given in the statement.
We have
\Bsist
  && - \intQt \zeta \phistar
  \leq \bigl( \norma{\zeta_1}_{\L2H} + \norma{\zeta_2}_{\L2H} \bigr) \norma\phistar_{\L2H}
  \non
  \\
  && \leq c \bigl( \norma{\phi_1-\phistar_1}_{\L2H} + \norma{\phi_2-\phistar_2}_{\L2H} + 1 \bigr) \norma\phistar_{\L2H} 
  \non
  \\
  && \leq c \, \norma\phistar_{\L2H} \,.
  \non
\Esist
Therefore, by also applying the Schwarz and Young inequalities, 
recalling the continuous embeddings $H\subset\Vp$ and $H\subset\Vzp$ 
(in~the two case of the boundary conditions) and rearranging, 
we deduce~that
\Bsist
  && \normaH{\phi(t)}^2
  + \intQt |\nabla\phi|^2
  \non
  \\
  && \leq c \, \normaH\phiz^2
  + c \, \norma g_{\L2H}^2
  + c \intQt |\phi|^2
  + c \, \norma\phistar_{\L2H} \,.
  \label{forcontdep}
\Esist
At this point, we can apply the Gronwall lemma and obtain the inequality \eqref{contdep}.\QED

\Brem
\label{RemD}
We come back to the case of the Dirichlet boundary conditions for the chemical potential
and make a remark concerning possibly different Dirichlet data.
For clarity, we denote by $\muG^i$, $\muH^i$ and $\mu^i$ the Dirichlet data, 
the corresponding harmonic extensions and the original chemical potentials associated to the given couple of data.
Then, the meaning of $\mu_i$ in the proof is $\mu_i:=\mu^i-\muG^i$
and the one of $g_i$ is the difference between the original $g_i$ and~$\muH^i$.
Hence, the \rhs\ of \eqref{contdep} contains both the norm in $\L2H$ of \gianni{the} difference of the original $g_1$ and $g_2$
and the one of $\muH^1$ and~$\muH^2$.
The latter is then estimated by the norm in $\L2{\HxG{1/2}}$ (or by a weaker norm) of $\muG^1-\muG^2$.
This means that $\phi$ continuously depends also on the Dirichlet datum for the chemical potential.
\Erem


\section{Approximation}
\label{APPROXIMATION}
\setcounter{equation}{0}

The method we use for the existence part of Theorem~\ref{Wellposedness}
consists in performing suitable a~priori estimates on the solution to an approximating problem
and using compactness and monotonicity arguments.
This section is devoted to the approximating problem.

We introduce the Moreau-Yosida regularizations $\betaeps$, $\Seps$, $\Betaeps$ and $\hatSeps$
of the graphs $\beta$ and $S$ and of their primitives $\Beta$ and~$\hatS$ at the level $\eps\in(0,1)$,
and replace $\beta$ and $S$ in \Pbl\ by $\betaeps$ and~$\Seps$, respectively.
We obtain the problem of finding $\soluzeps$ satisfying \accorpa{regphi}{regmu} as well~as
\Bsist
  && \iO \dt\phieps(t) \, v 
  + \iO \nabla\mueps(t) \cdot \nabla v
  = 0 
  \quad \hbox{\aat\ and every $v\in\calV$}
  \qquad
  \label{primaeps}
  \\
  \separa
  && \tau \iO \dt\phieps(t) \, v
  + \iO \nabla\phieps(t) \cdot \nabla v
  + \iO \bigl( \betaeps(\phieps(t)) + \pi(\phieps(t)) \bigr) v
  \non
  \\
  && \quad {}
  + \iO \Seps(\phieps(t)-\phistar(t)) v
  \non
  \\
  && = \iO \bigl( \mueps(t) + g(t) \bigr) v
  \quad \hbox{\aat\ and every $v\in V$}
  \label{secondaeps}
  \\
  \separa
  && \phieps(0) = \phiz 
  \label{cauchyeps}
\Esist
\Accorpa\Pbleps primaeps cauchyeps
where $\calV$ is either $V$ or $\Vz$ according to Notation~\ref{Notation}.

\Bthm
\label{Wellposednesseps}
Problem \Pbleps\ has a unique solution $\soluzeps$ 
satisfying the regularity requirements \accorpa{regphi}{regmu}.
\Ethm

\Bdim
Uniqueness for $\phieps$ is still given by Section~\ref{UNIQUENESS}
since \Pbleps\ is a particular case of problem \Pbl.
Indeed, $\betaeps$~and $\Seps$ satisfy the assumptions
we have required for $\beta$ and~$S$.
This is clear for~$\betaeps$.
As for~$\Seps$, we recall that $\normaH{\Seps(v)}\leq\normaH{\Sz(v)}$ and $\Sz(v)\in S(v)$ for every $v\in H$,
so that the linear growth condition~\eqref{hpS} yields 
\Beq
  \normaH{\Seps(v)} \leq C_S (\normaH v + 1)
  \quad \hbox{for every $v\in H$} .
  \label{disugSeps}
\Eeq
Once uniqueness for $\phieps$ is established,
uniqueness for $\mueps$ trivially follows by comparison in \eqref{secondaeps}
since both $\betaeps$ and $\Seps$ are single-valued.

As for existence, we can follow the line of \cite{CGMR3} only in the case of the Dirichlet boundary conditions for~$\mueps$.
However, we have to take care on the dependence of $\phistar$ on time
($\phistar$~is constant in time in~\cite{CGMR3}) and on the general shape of~$S$.
On~the contrary, the treatment of the Neumann boundary conditions for $\mueps$ needs a different argument.
This is developed in the next lines
by performing a discretization by means of a Faedo--Galerkin scheme
corresponding to a Hilbert basis of eigenfunctions.

\step 
The case of the Neumann boundary conditions: discretization

Let $\graffe{\lambda_j}_{j\geq1}$ and $\graffe{e_j}_{j\geq1}$
be the sequence of the eigenvalues and an orthonormal system of corresponding eigenfunctions
of the Neumann problem for the Laplace equation,~i.e.,
\Bsist
  && 0 = \lambda_1 < \lambda_2 \leq \lambda_3 \leq \dots
  \aand \lim_{j\to\infty} \lambda_j = + \infty
  \non
  \\
  && e_j\in\calV
  \aand
  \iO \nabla e_j \cdot \nabla e_j
  = \lambda_j \iO e_j v
  \quad \hbox{for every $v\in\calV$ and $j=1,2,\dots$}
  \non
  \\
  && \iO e_i e_j = \delta_{ij}
  \quad \hbox{for $i,j=1,2,\dots$}
  \aand 
  \hbox{$\graffe{e_j}_{j\geq1}$ is a complete system in $H$}.
  \non
\Esist
Notice that we can take $e_1=\mO^{-1/2}$.
We set
\Beq
  \Vn := \Span \graffe{e_1,\dots,e_n}
  \quad \hbox{for $n=1,2,\dots$}
  \aand
  \Vinfty := \bigcup_{n=1}^\infty \Vn
  \label{defVn}
\Eeq
and observe that $\Vinfty$ is dense in both $V$ and~$H$.
The discretized problem consists in finding 
$\phin\in\H1\Vn$ and $\mun\in\L2\Vn$ satisfying
\Bsist
  && \iO \dt\phin(t) \, v 
  + \iO \nabla\mun(t) \cdot \nabla v
  = 0 
  \quad \hbox{\aat\ and every $v\in\Vn$}
  \qquad
  \label{priman}
  \\
  \separa
  && \tau \iO \dt\phin(t) \, v
  + \iO \nabla\phin(t) \cdot \nabla v
  + \iO \bigl( \betaeps(\phin(t)) + \pi(\phin(t)) \bigr) v
  \non
  \\
  && \quad {}
  + \iO \Seps(\phin(t) - \phistar(t)) v
  \non
  \\
  && = \iO \bigl( \mun(t) + g(t) \bigr) v
  \quad \hbox{\aat\ and every $v\in\Vn$}
  \label{secondan}
  \\
  \separa
  && \iO \phin(0) v = \iO \phiz \, v
  \quad \hbox{for every $v\in\Vn$} \,.
  \label{cauchyn}
\Esist
\Accorpa\Pbln priman cauchyn
Thus, $\phin$ and $\mun$ have to be expanded as
\Bsist
  && \phin(t) = \somma j1n \phij(t) e_j
  \aand
  \mun(t) = \somma j1n \etaj(t) e_j
  \non
  \\
  && \quad \hbox{for some $\phij\in H^1(0,T)$ and $\etaj\in L^2(0,T)$, \ $j=1,\dots,n$},
  \non
\Esist
and problem \Pbln\ written in terms of the column vectors
$\badphi:={}^t(\badphi_1,\dots,\badphi_n)$ and $\eta:={}^t(\eta_1,\dots,\eta_n)$ takes the form
\Bsist
  && \badphi'(t) + A \eta(t) = 0
  \quad \aat
  \label{primad}
  \\
  && \tau \badphi'(t) 
  + A \badphi(t) 
  + \beps(\badphi(t))
  + \sigmaeps(t,\badphi(t))
  \non
  \\
  && = \eta(t) + \gn(t)
  \quad \aat
  \label{secondad}
\Esist
with an additial initial condition for~$\badphi$,
where $A:=(a_{ij})\in\erre^{n\times n}$,
$\beps:=(\beps_i):\erre^n\to\erre^n$, $\sigmaeps:=(\sigmaeps_i):(0,T)\times\erre^n\to\erre^n$ 
and $\gn:=(\gn_i):(0,T)\to\erre^n$ are given~by
\Bsist
  && a_{ij}
  := \iO \nabla e_j \cdot \nabla e_i \,, \quad
  \beps_i(y_1,\dots,y_n)
  := \iO (\betaeps+\pi) \Bigl( {\textstyle\somma j1n y_j e_j} \Bigr) e_i
  \non
  \\
  && \sigmaeps_i(t;y_1,\dots,y_n)
  := \iO \Seps \Bigl( {\textstyle\somma j1n y_j e_j} - \phistar(t) \Bigr) e_i
  \aand
  \gn_i(t) 
  := \iO g(t) e_i \,.
  \non
\Esist
Since $\betaeps:\erre\to\erre$ is \Lip\ continuous,
$\beps$ is \Lip\ continuous. 
On the other hand, $\Seps:H\to H$ is \Lip\ continuous too
and the function $(y_1,\dots,y_n)\mapsto\normaH{\somma j1n y_j e_j}$ is a norm in~$\erren$.
Hence, by denoting by $L_\eps$ the \Lip\ constant of~$\Seps$, we have for $i=1,\dots,n$ and every $y,z\in\erren$  
\Bsist
  && |\sigmaeps_i(t,y) - \sigmaeps_i(t,z)|
  = \Bigl|
    \iO \Seps \bigl( {\textstyle\somma j1n y_j e_j} - \phistar(t) \bigr) e_i 
    - \iO \Seps \bigl( {\textstyle\somma j1n z_j e_j} - \phistar(t) \bigr) e_i 
  \Bigr|
  \non
  \\
  && \leq \Bigl\|
    \Seps \bigl( {\textstyle\somma j1n y_j e_j} - \phistar(t) \bigr) 
    - \Seps \bigl( {\textstyle\somma j1n z_j e_j} - \phistar(t) \bigr) 
  \Bigr\|_H
  \non
  \\
  && \leq L_\eps \, \bigl\|
    {\textstyle\somma j1n y_j e_j} 
    - {\textstyle\somma j1n z_j e_j}
  \bigr\|_H
  \approx |y-z| 
  \non
\Esist
that is, $\sigmaeps$ is a Charath\'eodory function
that is \Lip\ continuous with respect to $y$ uniformly in~$t$.
Moreover, due to our assumptions \eqref{hpg} and \eqref{hpphistar} on $g$ and~$\phistar$,
we have that $t\mapsto\sigmaeps(t;y)$ and $\gn$ belong to $\L2{\erre^n}$, the former for every $y\in\erre^n$.
By replacing $\eta$ in \eqref{primad} by the expression obtained by solving \eqref{secondad} for~$\eta$,
we obtain a Cauchy problem for the single equation (where $I\in\erre^{n\times n}$ is the identity matrix)
\Beq
  (I + \tau A) \badphi'(t)
  + A \bigl(
  A \badphi(t) 
  + \beps(\badphi(t))
  + \sigmaeps(t,\badphi(t))
  - \gn(t)
  \bigr)
  = 0
  \non
\Eeq
which can be put in its normal form since $I+\tau A$ is positive definite.
Thus, we obtain a unique $\badphi\in\H1{\erre^n}$
and \eqref{secondad} provides the corresponding $\eta\in\L2{\erre^n}$.
This proves that the discrete problem \Pbln\ has a unique solution with the desired regularity.
At this point, our aim is letting $n$ tend to infinity.
To this end, we perform an a~priori estimate.
In the sequel, $c$~stands for different constants independent of~$n$
(and possibly depending on the structure of our system, $\Omega$, $T$, the data and~$\eps$, 
which is fixed in the whole proof of Theorem~\ref{Wellposednesseps}).

\step
Estimate on the discrete solution

We test \eqref{priman} and \eqref{secondan} by the $\Vn$-valued functions
$\mun$ and~$\dt\phin$, respectively, sum up and integrate over~$(0,t)$.
Moreover, we add the same quantity $\frac12\iO|\phin(t)|^2=\intQt\phin\dt\phin+\frac12\iO|\phin(0)|^2$
to both sides.
After an obvious cancellation and a rearrangement, we obtain
\Bsist
  && \intQt |\nabla\mun|^2
  + \tau \intQt |\dt\phin|^2
  + \normaV{\phin(t)}^2
  + \iO \Betaeps(\phin(t))
  \non
  \\
  && = \normaV{\phin(0)}^2
  + \iO \Betaeps(\phin(0))
  + \intQt \bigl(
    g
    + \phin
    - \pi(\phin) 
    - \Seps(\phin-\phistar) 
  \bigr) \dt\phin \,.
  \qquad
  \label{perstiman}
\Esist
The last integral can be treated by \pier{using} the \Lip\ continuity of $\pi$ and $\Seps$
and applying the Schwarz and Young inequalities as follows
\Bsist
  && \intQt \bigl(
    g
    + \phin
    - \pi(\phin) 
    - \Seps(\phin-\phistar) 
  \bigr) \dt\phin 
  \non
  \\
  && \leq \frac \tau 2 \intQt |\dt\phin|^2
  + c \intQt ( |\phin|^2 + 1)
  + c \, \norma g_{\L2H}^2 
  + c \, \norma\phistar_{\L2H}^2 \,.
  \non
\Esist
If we prove that $\phin(0)$ is bounded in~$V$,
even the second term on the \rhs\ of \eqref{perstiman} remains bounded
(since $\Betaeps(r)$ grows at most as~$r^2$)
and we can apply the Gronwall lemma to conclude~that
\Beq
  \norma\phin_{\H1H\cap\pier{{}\L\infty V}}
  + \norma{\nabla\mun}_{\L2H}
  \leq c \,.
  \label{semistiman}
\Eeq
\pier{From \eqref{cauchyn} it follows} that $\phin(0)$ is the $H$-projection of $\phiz$ onto~$\Vn$, i.e.,
\Beq
  \phin(0) = \somma j1n \alpha_j e_j
  \quad \hbox{where the sequence $\graffe{\alpha_j}$ satisfies} \quad
  \phiz = \somma j1\infty \alpha_j e_j \,.
  \non
\Eeq
Now, we observe that~$W$ (see \eqref{defspazi}) can be characterized~as 
\Beq
  W = \graffe{v\in H:\ -\Delta v\in H, \ \dn v=0}
  = \Bigl\{
    {\textstyle \somma j1\infty c_j e_j : \ \somma j1\infty \lambda_j^2 |c_j|^2 < +\infty}
  \Bigr\}.
  \non
\Eeq
Since $V$ is the interpolation space $(W,H)_{1/2}$,
we also have that
\Bsist
  && V = \Bigl\{
    {\textstyle \somma j1\infty c_j e_j : \ \somma j1\infty \lambda_j |c_j|^2 < +\infty}
  \Bigr\}
  \quad \hbox{with equivalence of Hilbert norms}
  \quad
  \non
  \\
  && \normaV v
  \approx \norma v_\lambda
  := \Bigl(
    \normaH v^2 + {\textstyle \somma j1\infty \lambda_j |c_j|^2}
  \Bigr)^{1/2}
  \quad \hbox{if $v=\somma j1\infty c_je_j$}.
  \non
\Esist
Since $\phiz\in V$, we deduce that
\Beq
  \norma{\phin(0)}_\lambda^2 
  = \normaH{\phin(0)}^2 + {\textstyle \somma j1n \lambda_j |\alpha_j|^2}
  \leq \normaH \phiz^2 + {\textstyle \somma j1\infty \lambda_j |\alpha_j|^2} 
  = \norma\phiz_\lambda^2 \,.
  \label{stimaphinz}
\Eeq
This concludes the proof of~\eqref{semistiman}.
We immediately infer that $\betaeps(\phin)$, $\pi(\phin)$ and $\Seps(\phin-\phistar)$ are bounded in $\L2H$.
At this point, we can test \eqref{priman} and \eqref{secondan} by the $\Vn$-valued functions
$\phin$ and~$-\mun$, respectively, sum up, integrate over~$(0,t)$
and notice that a cancellation occurs.
Hence, some rearrangement, the use of the previous estimates and the Young inequality lead~to
\Bsist
  && \frac 12 \normaH{\phin(t)}^2
  + \intQt |\mun|^2
  = \frac 12 \normaH{\phin(0)}^2
  \non
  \\
  && \quad {}
  + \intQt \bigl(
    \tau \dt\phin
    + (\betaeps+\pi)(\phin) 
    + \Seps(\phin-\phistar)
    - g
  \bigr) \mun 
  \non
  \\
  && \leq c + \frac 12 \intQt |\mun|^2 
  \non
\Esist
whence a bound for $\mun$ in $\L2H$.
Therefore, \eqref{semistiman}~is improved.
Namely, we have
\Beq
  \norma\phin_{\H1H\cap\pier{{}\L\infty V}}
  + \norma\mun_{\L2V}
  \leq c \,.
  \label{stiman}
\Eeq

\step
Conclusion

From \eqref{stiman}, the Aubin--Lions lemma
(see, e.g., \cite[Thm.~5.1, p.~58]{Lions})
and the \Lip\ continuity of $\betaeps$, $\pi$ and~$\Seps$,
we deduce that
\Bsist
  && \phin \to \phieps
  \quad \hbox{weakly \pier{star} in $\H1H\cap\pier{{}\L\infty V}$}
  \label{convphin}
  \\
  && \mun \to \mueps
  \quad \hbox{weakly in $\L2V$}
  \label{convmun}
  \\
  && \phin \to \phieps , \quad
  (\betaeps+\pi)(\phin) \to (\betaeps+\pi)(\phieps)
  \non
  \\
  && \aand
  \Seps(\phin-\phistar) \to \Seps(\phieps-\phistar)
  \quad \hbox{strongly in $\L2H$} 
  \label{nonlinn}
\Esist
as $n$ tends to infinity, at least for~a (not relabeled) subsequence.
\pier{Actually, $\phin$ converges to $\phieps$ strongly in $\C0H$ due to 
\eqref{convphin} and the generalized Ascoli theorem 
(see, e.g., \cite[Sect.~8, Cor.~4]{Simon}).} 
We show that $\soluzeps$ is the solution to \Pbleps\ we are looking for.
First, we have that $\phieps(0)=\phiz$
since $\phin(0)$ converges \pier{both to $\phieps(0)$ 
and to $\phiz$ strongly in~$H$, the latter since $\Vinfty$ 
(see~\eqref{defVn}) is dense in~$H$ (in~fact, 
\eqref{stimaphinz} ensures that $\phin(0)\to \phiz$ strongly in $V$).}
Now, fix $m\geq1$, take any $v\in\L2\Vm$ and assume that $n\geq m$.
Then, $\Vm\subset\Vn$ so that both \eqref{priman} and \eqref{secondan}
can be tested by $v(t)$ and then integrated over~$(0,T)$.
We thus obtain
\Bsist
  && \intQ \dt\phin \, v 
  + \intQ \nabla\mun \cdot \nabla v
  = 0 
  \non
  \\
  && \tau \intQ \dt\phin \, v
  + \intQ \nabla\phin \cdot \nabla v
  + \intQ \bigl( \betaeps(\phin) + \pi(\phin) \bigr) v
  + \intQ \Seps(\phin - \phistar) v
  \non
  \\
  && = \intQ \bigl( \mun + g \bigr) v 
  \non
\Esist
for arbitrary choices of $v\in\L2\Vm$ in the above equations
and every $n\geq m$.
By letting $n$ tend to infinity, we deduce that
\Bsist
  && \intQ \dt\phieps \, v 
  + \intQ \nabla\mueps \cdot \nabla v
  = 0 
  \label{intprimaeps}
  \\
  && \tau \intQ \dt\phieps \, v
  + \intQ \nabla\phieps \cdot \nabla v
  + \intQ \bigl( \betaeps(\phieps) + \pi(\phieps) \bigr) v
  + \intQ \Seps(\phieps - \phistar) v
  \non
  \\
  && = \intQ \bigl( \mueps + g \bigr) v .
  \label{intsecondaeps}
\Esist
In both equations, $v\in\L2\Vm$ is arbitrary
and $m$ is arbitrary too.
Take now any $\Vinfty$-valued step function.
Then, $v\in\L2\Vm$ for some~$m$ 
so that \eqref{intprimaeps} and \eqref{intsecondaeps} hold for~$v$.
Since the set of such functions is dense in~$\L2V$,
we conclude that both \eqref{intprimaeps} and \eqref{intsecondaeps} hold for every $v\in\L2V$.
But this is equivalent to say that \eqref{primaeps} and \eqref{secondaeps} hold true.
This concludes the proof in the case of the Neumann conditions.

\step
The case of the Dirichlet boundary conditions

In this case, the argument used above does not work, 
since we should use two different systems of eigenfunctions
corresponding to the different boundary conditions for $\phieps$ and~$\mueps$.
Therefore, we follow the ideas of~\cite{CGMR3} and eliminate~$\mueps$.
By recalling the definition \eqref{defD} of~$\calD$, we write \eqref{primaeps} in the equivalent form
\Beq
  \mueps = - \calD(\dt\phieps) 
  \label{muepsD}
\Eeq
and replace $\mueps$ in \eqref{secondaeps} by this expression.
We obtain
\Bsist
  && \iO \bigl( \tau \dt\phieps(t) + \calD(\dt\phieps(t)) \bigr) v 
  + \iO \nabla\phieps(t) \cdot \nabla v
  \non
  \\
  && \quad {}
  + \iO \bigl( \betaeps(\phieps(t)) + \pi(\phieps(t)) + \Seps(\phieps(t)-\phistar(t)) \bigr) v
  \non
  \\
  && = \iO g(t) v
  \quad \hbox{\aat\ and every $v\in V$} .
  \label{secondaepsD}
\Esist
We present this variational equation as an abstract nonlinear equation of parabolic type.
To this end, we introduce the bilinear form $\(\cdot,\cdot\)$ in $H\times H$ by setting
\Bsist
  && \(u,v\)
  := \iO ( Bu ) v 
  \quad \hbox{for $u,v\in H$},
  \quad \hbox{where} 
  \non
  \\
  && B := \tau I + \calD \,,
  \quad \hbox{$I$ being the identity map of $H$}
  \non
\Esist
and observe that, due to the definition \eqref{defD} of~$\calD$, 
this form is continuous and symmetric and satisfies for every $v\in H$
\Beq
  \(v,v\)
  = \tau \, \normaH v^2 
  + \iO (\calD v) v
  = \tau \, \normaH v^2 
  + \iO |\nabla\calD v|^2
  \geq \tau \, \normaH v^2 \,.
  \non
\Eeq
Hence, it is an equivalent inner product in~$H$.
This also shows that the operator $B$ is an isomorphism from $H$ into itself.
Therefore, the variational equation \eqref{secondaepsD} can be written~as
\Bsist
  && \( \dt\phieps,v \)
  + \iO \nabla\phieps\cdot \nabla v
  + \( B^{-1} \bigl( \betaeps(\phieps) + \pi(\phieps) + \Seps(\phieps-\phistar) \bigr) , v \)
  \non
  \\
  && = \( B^{-1}g , v \)
  \quad \hbox{\aet, for every $v\in V$} 
  \non
\Esist
that is, as the abstract equation in $\Vp$
\Beq
  \frac d{dt} \, \phieps(t) + A \phieps(t) + \Feps(t,\phieps(t)) = B^{-1} g(t) 
  \quad \aat
  \label{abstract}
\Eeq
provided that the Hilbert triplet $(V,H,\Vp)$ is constructed 
starting from the new inner product $\(\cdot,\cdot\)$ rather than the standard one
(i.e.,~in order that $\<u,v>=\(u,v\)$ for every $u\in H$ and $v\in V$)
and the continuous linear operator $A:V\to\Vp$ and the function $\Feps:(0,T)\times H\to H\subset\Vp$
are defined as follows
\Bsist
  && \< Au, v > := \iO \nabla u \cdot \nabla v
  \quad \hbox{for $u,v\in V$}
  \non
  \\
  && \Feps(t,v)
  := \pier{ B^{-1} }\bigl( \betaeps(v) + \pi(v) + \Seps(v-\phistar(t)) \bigr)
  \quad \hbox{\aat\ and $v\in H$}.
  \non
\Esist
So, it is sufficient to solve the Cauchy problem for \eqref{abstract}
associated \pier{with} the initial condition \eqref{cauchyeps} and then recover $\mueps$ by means of~\eqref{muepsD}.
The Cauchy problem just mentioned has a unique solution $\phieps\in\H1H\cap\pier{{}\L\infty V}$ by the following reasons.
First, $A$~is linear and continuous and $A+I$ (where $I:V\to\Vp$ is the embedding) is~coercive.
Next, since $\phistar\in\L2H$ and $\betaeps$, $\pi$ and $\Seps$ are \Lip\ continuous,
$\Feps$~is a Charat\'eodory function that is \Lip\ continuous with respect to~$v$ uniformly in~$t$.
Finally, the \rhs\ belongs to $\L2H$ and $\phiz\in V$.
For a more detailed proof, one could discretize \eqref{abstract} with a Faedo--Galerkin scheme,
as we did before.
Indeed, as $\mueps$ has been eliminated, one can use just one system of eigenfunction,
namely, the same we have introduced to solve the problem in the case of the Neumann boundary conditions.
This concludes the proof.
\Edim


\section{Existence for the generalized problem}
\label{EXISTENCE}
\setcounter{equation}{0}

We start from the approximating problem \Pbleps\ and perform some a priori estimates on its solution $\soluzeps$.
In the whole section, we use the same symbol $c$ to denote constants
that only \pier{depend} on the structure of the problem, the data, $\Omega$ and~$T$,
and can possibly be different from each other (even in the same chain of equalities or inequalities).
We~stress that the values of $c$ do not depend on~$\eps$.

\step
First a priori estimate

Assume first that $\calV=V$.
Then, by taking $v=\mO^{-1}$ in~\eqref{primaeps}, we see that
$\mean(\dt\phieps)=0$ so that $\calN(\dt\phieps)$ is well defined.
Moreover, it belongs to $\L2V$.
If instead $\calV=\Vz$ one can consider $\calD(\dt\phieps)$,
which is well defined and belongs to~$\L2\Vz$.
Hence, in both cases, we can test \eqref{primaeps} written at the time $s$ by $\mueps(s)+\calL(\dt\phieps(s))$
where $\calL=\calN$ if $\calV=V$ and $\calL=\calD$ if $\calV=\Vz$.
Then, we integrate the resulting equality over $(0,t)$ with respect to~$s$.
Similarly, we test \eqref{secondaeps} by $2\dt\phieps$ and integrate in time.
Then, we sum the equalities obtained this way to each other
and notice that two cancellations occur:
one of them is obvious and the other is due to the definitions of $\calN$ and $\calD$ given by \eqref{defN} and~\eqref{defD}.
Finally, we add the same quantity 
$\iO|\phieps(t)|^2=\iO|\phiz|^2+2\intQt\phieps\dt\phieps$ to both sides for convenience.
Hence, by recalling~\eqref{normaVp} and \eqref{normaVzp}, we obtain
\Bsist
  && \iot \normaVp{\dt\phieps(s)}^2 \ds
  + \intQt |\nabla\mueps|^2
  \non
  \\
  && \quad {}
  + 2 \tau \intQt |\dt\phieps|^2
  + \iO |\nabla\phieps(t)|^2
  + 2 \iO \Betaeps(\phieps(t))
  + \iO |\phieps(t)|^2
  \non
  \\
  \separa
  && = \iO |\nabla\phiz|^2
  + 2 \iO \Betaeps(\phiz)
  + \iO |\phiz|^2
  + 2 \intQt \bigl(
    g
    - \pi(\phieps)
    \pier{{}+ \phieps}
  \bigr) \dt\phieps 
  \non
  \\
  && \quad {}
  - 2 \intQt \Seps(\phieps-\phistar) \dt\phieps \,.
  \label{perprimastima}
\Esist
All the terms on the \lhs\ are nonnegative.
As for the first line on the \rhs, 
we account for the well-known inequality $\Betaeps(r)\leq\Beta(r)$
(which holds for every $r\in\erre$),
the linear growth of~$\pi$ (see~\eqref{hpPi})
and the Schwarz and Young inequalities.
Hence, we can estimate it from above by the following expression
\Bsist
  && \normaV\phiz^2
  + 2 \iO \Beta(\phiz)
  + \tau \intQt |\dt\phieps|^2
  + c \intQt \bigl(
    |g|^2 + |\phieps|^2 + 1
  \bigr).
  \non
\Esist
Next, we use the Young inequality, the uniform linear growth condition \eqref{disugSeps}
and the assumption \eqref{hpphistar} on~$\phistar$.
We obtain
\Bsist
  && - \intQt \Seps(\phieps-\phistar) \dt\phieps
  \leq \frac \tau 2 \intQt |\dt\phieps|^2
  + c \intQt (|\phieps|^2 + |\phistar|^2) 
  \non
  \\
  && \leq \frac \tau 2 \intQt |\dt\phieps|^2
  + c \intQt |\phieps|^2 
  + c \,.
  \label{termineSeps}
\Esist
At this point, we come back to \eqref{perprimastima}, combine with these estimates,
rearrange and apply the Gronwall lemma.
We conclude~that
\Beq
  \norma\phieps_{\H1H\cap\L\infty V} 
  + \norma{\nabla\mueps}_{\L2H}
  \leq c \,.
  \label{primastima}
\Eeq
\pier{Hence, by accounting also for the \Lip\ continuity of~$\pi$, the inequality \eqref{disugSeps}}
and the $L^2$ summability of~$\phistar$, we deduce~that
\Beq
  \norma\phieps_{\L\infty H}
  + \norma{\pi(\phieps)}_{\L\infty H}
  + \norma{\Seps(\phieps-\phistar)}_{\L2H} 
  \leq c \,.
  \label{daprimastima}
\Eeq

\step
Second a priori estimate

Our aim is to improve a part of \eqref{primastima} by showing that
\Beq
  \norma\mueps_{\L2\calV}
  \leq c \,.
  \label{stimamueps}
\Eeq
If $\calV=\Vz$ this trivially follows from \eqref{primastima} and the second Poincar\'e inequality~\eqref{poincare}.
If instead $\calV=V$, we can deduce \eqref{stimamueps} from the first Poincar\'e inequality 
provided we have an estimate of the mean value $\mean\mueps$ in~$L^2(0,T)$.
Hence, we prepare \pier{a} pointwise estimate in this direction.

Since $\calV=V$, one of our assumption is~\eqref{hpmz}.
The following argument (which owes to \cite[Appendix, Prop.~A.1]{MiZe},
see also \cite[p.~908]{GiMiSchi} for a detailed proof) is~used in several papers.
We repeat it here for the reader's convenience. 
By~\eqref{hpmz}, we have for some $\delta_0>0$ depending only on $\beta$ and~$\mz$
\Beq
  \betaeps(s) (r-\mz)
  \geq \delta_0 |\betaeps(r)| - \delta_0^{-1}
  \quad \hbox{for every $r\in\erre$ and every $\eps\in(0,1)$}.
  \label{trickMZ}
\Eeq
Now, we test \eqref{secondaeps} by~$\phieps-\mz$ and avoid time integration.
Thus, we argue for $t$ fixed.
However, for simplicity, we do not write the time $t$ for a while.
We have \aet
\Bsist
  && \pier{\delta_0}\iO |\betaeps(\phieps)| - \delta_0^{-1} \, \mO
  \leq \iO \nabla\phieps \cdot \nabla(\phieps-\mz)
  + \iO \betaeps(\phieps) (\phieps-\mz)
  \non
  \\
  && = \iO \mueps (\phieps-\mz)
  + \iO \bigl(
    g
    - \tau \dt\phieps
    - \pi(\phieps)
    - \Seps(\phieps-\phistar)
  \bigr) (\phieps-\mz) .
  \label{perpt}
\Esist
We recall that $\mean(\phieps-\mz)=0$ \aet\pier{, thus we can take advantage 
of that and} apply the first Poincar\'e inequality~\eqref{poincare} \pier{to 
$\mueps-\mean\mueps$. In fact, using \eqref{daprimastima} as well,} we have
\Bsist
  && \iO \mueps (\phieps-\mz)
  = \iO (\mueps-\mean\mueps) (\phieps-\mz)
  \leq \normaH{\mueps-\mean\mueps} \, \normaH{\phieps-\mz}  
  \non
  \\
  && \leq c \, \normaH{\nabla\mueps} \, \normaH{\phieps-\mz} 
  \leq c \, \normaH{\nabla\mueps}  \,.
  \non
\Esist
As for the rest of the \rhs\ of \eqref{perpt},
we use the Schwarz inequality and the estimates in $\L\infty H$ available from~\eqref{daprimastima}.
Hence, we obtain
\Bsist
  && \norma{\betaeps(\phieps)}_1
  \leq c \, \bigl(
    \normaH{\nabla\mueps} 
    + \normaH g
    + \normaH{\dt\phieps}
    + \normaH{\Seps(\phieps-\phistar)}
    + 1 
  \bigr) 
  \non
\Esist
\aet.
By taking $v=\mO^{-1}$ in \eqref{secondaeps} as before, using the inequality just obtained
and estimating the other $L^1$-norms by the corresponding $H$-norms, 
we deduce that
\Beq
  |\mean\mueps|
  \leq c \, \bigl(
    \normaH{\nabla\mueps} 
    + \normaH g
    + \normaH{\dt\phieps}
    + \normaH{\Seps(\phieps-\phistar)} 
    + 1
  \bigr) 
  \label{pointwise}
\Eeq
\aet.
This is the desired pointwise estimate.

By squaring \eqref{pointwise}, integrating over~$(0,T)$ and taking the square root,
we deduce an inequality.
The \lhs\ is the norm of $\mean\mueps$ in $L^2(0,T)$
and the \rhs\ is what one obtains by replacing 
every $H$-norm in the \rhs\ of \eqref{pointwise} by the norm in~$\L2H$.
Hence, by using \eqref{primastima} and~\eqref{daprimastima},
we deduce that
\Beq
  \norma{\mean\mueps}_{L^2(0,T)} \leq c \,.
  \non
\Eeq
Hence, \eqref{stimamueps} follows.

\step
Third a priori estimate

We test \eqref{secondaeps} by $\betaeps(\phieps)$ and integrate over~$(0,T)$.
We obtain
\Bsist
  && \intQ \betaeps'(\phieps) |\nabla\phieps|^2
  + \intQ |\betaeps(\phieps)|^2
  \non
  \\
  && = \intQ \bigl(
    \mueps + g - \tau \dt\phieps - \pi(\phieps) - \Seps(\phieps-\phistar)
  \bigr) \betaeps(\phieps).
  \label{perterzastima}
\Esist
From the previous estimates we immediately conclude that
\Beq
  \norma{\betaeps(\phieps)}_{\L2H} \leq c \,.
  \label{terzastima}
\Eeq

\step
Conclusion

By accounting for \accorpa{primastima}{stimamueps} and~\eqref{terzastima}
and owing to well-known weak compactness results,
we have in both cases $\calV=V$ and $\calV=\Vz$ (along a subsequence)
\Bsist
  && \phieps \to \phi
  \quad \hbox{weakly in $\H1H\cap\L2{\pier{V}}$}
  \label{convphieps}
  \\
  && \mueps \to \mu
  \quad \hbox{weakly in $\L2\calV$}
  \label{convmueps}
  \\
  && \betaeps(\phieps) \to \xi
  \quad \hbox{weakly in $\L2H$}
  \label{convxieps}
  \\
  && \Seps(\phieps-\phistar) \to \zeta
  \quad \hbox{weakly in $\L2H$} 
  \label{convzetaeps}   
\Esist
for some quadruplet $\soluz$ with the regularity specified by the convergence properties.
We~claim that $\soluz$ is the solution to problem \Pbl\ we are looking for.
By letting $\eps$ tend to zero in~\pier{\eqref{primaeps}}, we clearly find that
\Beq
  \intQ \dt\phi \, v
  + \intQ \nabla\mu \cdot \nabla v
  = 0 
  \non
\Eeq
for every $v\in\L2\calV$, and this is equivalent to~\eqref{prima}.
Moreover, by \eqref{convphieps}, $\phieps$~converges to $\phi$ weakly in $\C0H$, 
whence we deduce that $\phieps(0)$ converge to $\phi(0)$ weakly in~$H$
and conclude that $\phi(0)=\phiz$.
Furthermore, we observe that the already mentioned Aubin--Lions lemma implies that
\Beq
  \phieps \to \phi
  \quad \hbox{strongly in $\L2H$}.
  \label{strongphieps}
\Eeq
By \Lip\ continuity, $\pi(\phieps)$ converges to $\pi(\phi)$ in the same topology.
It follows that
\Beq
  \tau \intQ \dt\phi \, v
  + \intQ \nabla\phi \cdot \nabla v
  + \intQ \bigl( \xi + \pi(\phieps) + \zeta \bigr) v
  = \intQ \bigl( \mueps + g \bigr) v 
  \non
\Eeq
for every $v\in\L2\calV$, and this is equivalent to \eqref{seconda}.
Finally, by accounting for the strong convergence \eqref{strongphieps}
(which trivially implies the strong convergence of $\phieps-\phistar$ to $\phi-\phistar$ in the same topology)
and the weak convergence given by \accorpa{convxieps}{convzetaeps},
we can apply, e.g., \cite[Lemma~2.3, p.~38]{Barbu}
and conclude that $\xi\in\beta(\phi)$ and $\zeta\in S(\phi-\phistar)$ \aeQ.
Hence, the proof is complete.


\section{Existence of a sliding mode}
\label{SLIDING}
\setcounter{equation}{0}

This section is devoted to the proof of Theorem~\ref{Sliding}.
Hence, $S$~is the graph induced in $H$ by $\rho\sign$ (see~\eqref{IdefS} and \accorpa{hphatS}{defsign})
and the data satisfy the further assumptions \HPdatibis\ as well.
Moreover, we choose $\calV=\Vz$ once and for all
in both the original problem \Pbl\ and the approximating problem \Pbleps,
since we only treat the Dirichlet boundary conditions for the chemical potential.
The proof of our result relies on a number of careful a~priori estimates 
on the solution $\soluzeps$ to the approximating problem
and a comparison argument involving the solution of an ODE.
In performing the a priori estimates, we have to take a particular care on 
the dependence of the constants on~$\rho$.
Thus, from now on, the symbol $c$ denotes (possibly different) constants 
that do not depend on $\rho$ and~$\eps$.
We use capital letters (mainly with indices like~$C_i$) to~denote some constants we want to refer~to.
The quantities that such constants depend on are specified in introducing them.
Our method is close to the ideas of~\cite{CGMR3} but it is different
(in~particular, we have to take care on the dependence of $\phistar$ on time).
For this reason, we present the whole argument with all the details that are necessary.

First of all, \pier{recalling} the properties  \eqref{hpphistarbis} on~$\phistar$,
we \pier{can} rewrite equations \accorpa{primaeps}{secondaeps} 
in \pier{terms} of both $\phieps$ and the auxiliary unknown
\Beq
  \chieps := \phieps - \phistar \,.
  \label{defchieps}
\Eeq
The new equations read
\Bsist
  && \iO \dt\chieps \, v
  + \iO \nabla\mueps \cdot \nabla v
  = - \iO \dt\phistar \, v
  \quad \hbox{for every $v\in\Vz$}
  \label{primaepsbis}
  \\
  && \tau \iO \dt\chieps \, v
  + \iO \nabla\chieps \cdot \nabla v
  + \iO \bigl(
    \betaeps(\phieps) + \pi(\phieps) + \Seps(\chieps)
  \bigr) v
  \non
  \\
  && = \iO \bigl( \mueps + g - \tau \dt\phistar \bigr) v
  - \iO \nabla\phistar \cdot \nabla v
  \quad \hbox{for every $v\in V$}
  \label{secondaepsbis}
\Esist
both \aet, and the initial condition for $\chieps$ is given by
\Beq
  \chieps(0)
  = \chiz := \phiz - \phistar(0) .
  \label{cauchyepsbis}
\Eeq
In equation \eqref{secondaepsbis}, $\Seps:H\to H$ is the map associated to $\seps:=\rho\signeps$,
that~is\pier{,}
\Beq
  \hbox{for $v\in H$}, \quad \hbox{$\Seps(v)$ is the function} \quad
  x \mapsto \seps(v(x)) = \rho\signeps v(x), \quad x\in\Omega ,
  \non
\Eeq
where $\signeps$ is the Yosida approximation of~$\sign$.
We \pier{point out} that
\Beq
  \signeps(0) = 0 
  \aand
  |\signeps(r)| \leq 1
  \quad \hbox{for every $r\in\erre$} .
  \label{propsigneps}
\Eeq
We also set for convenience 
\Beq
  \hatseps(r) := \int_0^r \seps(r') \, dr'
  \quad \hbox{for $r\in\erre$}
  \aand
  \hatSeps(v) := \iO \hatseps(v)
  \quad \hbox{for $v\in H$}
  \label{defhatSeps}
\Eeq
and notice that $\hatSeps$ actually is the primitive of~$\Seps$.
We also remark that
\Beq
  0 \leq \hatseps(r) \leq \rho |r|
  \aand
  0 \leq \hatSeps(v) \leq \hatS(v) = \rho \iO |v|
  \label{prophatSeps}
\Eeq
for every $r\in\erre$ and every $v\in H$, respectively.

\step
Fourth a priori estimate

We test \eqref{primaepsbis} and \eqref{secondaepsbis} written at the time~$s$ 
by $\calD(\chieps(s))$ and~$\chieps(s)$, respectively.
Then, we sum up, integrate over $(0,t)$ with respect to~$s$ 
and notice a cancellation due to the definition of~$\calD$ in the terms involving~$\mueps$.
By recalling \eqref{propND}, we obtain
\Bsist
  && \frac 12 \, \normaVp{\chieps(t)}^2
  + \frac \tau 2 \iO |\chieps(t)|^2
  + \intQt |\nabla\chieps|^2
  + \intQt \betaeps(\phieps) (\phieps-\phistar)
  + \rho \intQt \signeps(\chieps) \chieps
  \non
  \\
  && = - \intQt \dt\phistar \, \calD(\chieps)
  + \intQt \bigl( g - \pi(\chieps+\phistar) - \tau \dt\phistar \bigl) \chieps
  - \intQt \nabla\phistar \cdot \nabla\chieps \,.
  \non
\Esist
The last term on the \lhs\ is nonnegative since $\signeps$ is monotone and vanishes on the origin.
For the one involving $\betaeps$ we use the convexity inequality for $\Betaeps$ as follows
\Beq
  \betaeps(\phieps) (\phieps-\phistar)
  \geq \Betaeps(\phieps) - \Betaeps(\phistar)
  \non
\Eeq
and notice that
\Beq
  \iO \Betaeps(\phieps)
  \geq 0 
  \aand
  \intQt \Betaeps(\phistar)
  \leq \intQt \Beta(\phistar)
  \leq \intQ \Beta(\phistar)
  = c
  \non
\Eeq
where $c$ actually is finite as a consequence \pier{of the assumptions~\eqref{hpphistarbis}: indeed, it suffices that 
$\phistar, \betaz (\phistar) \in L^2(Q)$.}
Finally, all the terms on \rhs\ can be easily treated with the Young inequality
since $\normaH{\calD v}\leq c\,\normaH v$,
as a trivial consequence of the continuity of $\calD$ from $\Vzp$ into~$\Vz$.
Hence, by applying the Gronwall lemma, we~conclude~that
\Beq
  \norma\chieps_{\L\infty H\cap\L2V}
  \leq c \,.
  \label{stimabase}
\Eeq
By recalling \eqref{hpphistarbis} on $\phistar$ and the \Lip\ continuity of~$\pi$, we deduce that
\Beq
  \norma\phieps_{\L\infty H\cap\L2V}
  \leq c
  \aand
  \norma{\pi(\phieps)}_{\L\infty H\cap\L2V}
  \leq c \,.
  \label{dastimabase}
\Eeq
In the above estimates, $c$ does not depend on~$\rho$, according to our general rule.

Now, we revisit the derivation of one of the estimates of Section~\ref{EXISTENCE}
in order to explicit the dependence of the constant on $\rho$ under the new assumption on~$S$.

\step
First a priori estimate revisited

We come back to \eqref{perprimastima} and replace the treatment \pier{as in} \eqref{termineSeps} by a different argument.
Owing to \eqref{hpphistarbis} and \accorpa{propsigneps}{prophatSeps}, we can write
\Bsist
  && - 2 \intQt \Seps(\phieps-\phistar) \dt\phieps
  = - 2 \intQt \seps(\chieps) (\dt\chieps + \dt\phistar)
  \non
  \\
  && = - 2 \iO \hatseps(\chieps(t)) 
  + 2 \iO \hatseps(\chiz)
  - 2 \intQt \seps(\chieps) \dt\phistar 
  \non
  \\
  && \leq 2 \rho \iO |\chiz| 
  + 2 \rho \intQ |\dt\phistar|
  \leq c \, \rho \,.
  \non
\Esist
Hence, by recalling the derivation of \eqref{primastima} and \pier{in view of} \eqref{dastimabase}, \pier{we have now the estimate}
\Beq
  \pier{\norma{\phieps}_{\L\infty V}+ {}}\norma{\dt\phieps}_{\L2H}
  + \norma{\nabla\mueps}_{\L2H}
  \leq c \, (1+\rho^{1/2}) .
  \label{primastimanuova}
\Eeq
\pier{Thanks} to the Poincar\'e inequality, we deduce~that
\Beq
  \norma\mueps_{\L2\Vz}
  \leq c \, (1+\rho^{1/2}) \,.
  \label{daprimastimanuova}
\Eeq
In particular, $\mueps$ is bounded in $\L2H$ by the same constant.

\step
Fifth a priori estimate

We proceed formally for brevity
(for the correct argument, one could follow the idea suggested at the end of the proof of Theorem~\ref{Wellposednesseps},
i.e., the discretization of the abstract equation \eqref{abstract} by a Faedo--Galerkin scheme,
and then \pier{performing} the proper estimates on the discrete solution).
So, we formally differentiate \accorpa{primaeps}{secondaeps} with respect to time 
and obtain the equations written below,
which hold \aet\ and for every $v\in\Vz$ and every $v\in V$, respectively.
\Bsist
  && \iO \dt^2\phieps \, v 
  + \iO \nabla\dt\mueps \cdot \nabla v
  = 0 
  \label{dtprimaeps}
  \\
  && \tau \iO \dt^2\phieps \, v
  + \iO \nabla\dt\phieps \cdot \nabla v
  + \iO \betaeps'(\phieps) \dt\phieps \, v
  + \iO \dt(\seps(\chieps)) \, v
  \non
  \\
  && = \iO \bigl(
    \dt\mueps + \dt g - \pi'(\phieps) \dt\phieps
  \bigr) v \,.
  \label{dtsecondaeps}
\Esist
Then, we test the above equations written at the time $s$ 
by $\calD(\dt\phieps(s))$ and $\dt\phieps(s)$, respectively,
we sum up and integrate with respect to $s$ over~$(0,t)$.
Due to the definition \eqref{defD} of~$\calD$, a cancellation occurs.
Hence, by recalling~\eqref{propND}, we have
\Bsist
  && \frac 12 \, \normaVp{\dt\phieps(t)}^2
  + \frac \tau 2 \iO |\dt\phieps(t)|^2
  + \intQt |\nabla\dt\phieps|^2
  \non
  \\
  && \quad {}
  + \intQt \betaeps'(\phieps) |\dt\phieps|^2
  + \intQt \dt(\seps(\chieps)) \dt\phieps
  \non
  \\
  && = \frac 12 \, \normaVp{\dt\phieps(0)}^2
  + \frac \tau 2 \iO |\dt\phieps(0)|^2
  + \intQt \bigl(
    \dt g - \pi'(\phieps) \dt\phieps
  \bigr) \dt\phieps \,.
  \label{testquintastima}
\Esist
The term involving $\seps$ is treated this way
\Bsist
  && \intQt \dt(\seps(\chieps)) \, \dt\phieps
  = \intQt \dt(\seps(\chieps)) \, \dt\chieps
  + \intQt \dt(\seps(\chieps)) \, \dt\phistar
  \non
  \\
  && = \intQt \seps'(\chieps) |\dt\chieps|^2
  - \intQt \seps(\chieps) \, \dt^2\phistar
  + \iO \seps(\chieps(t)) \dt\phistar(t)
  - \iO \seps(\chiz) \dt\phistar(0) \,.
  \non
\Esist
Hence, by substituting in \eqref{testquintastima} 
and then ignoring the nonnegative terms containing $\betaeps'$ and~$\seps'$, 
we obtain
\Bsist
  && \frac 12 \, \normaVp{\dt\phieps(t)}^2
  + \frac \tau 2 \iO |\dt\phieps(t)|^2
  + \intQt |\nabla\dt\phieps|^2
  \non
  \\
  && \leq \frac 12 \, \normaVp{\dt\phieps(0)}^2
  + \frac \tau 2 \iO |\dt\phieps(0)|^2
  + \intQt \bigl(
    \dt g - \pi'(\phieps) \dt\phieps
  \bigr) \dt\phieps 
  \non
  \\
  && \quad {}
  + \intQt \seps(\chieps) \, \dt^2\phistar
  - \iO \seps(\chieps(t)) \dt\phistar(t)
  + \iO \seps(\chiz) \dt\phistar(0) \,.
  \qquad
  \label{perquintastima}
\Esist
We just have to estimate the \rhs.
The initial value $\dt\phieps(0)$ is read 
in equations \accorpa{primaeps}{secondaeps} written at $t=0$,~i.e.,
\Bsist
  && \iO \dt\phieps(0) \, v 
  + \iO \nabla\mueps(0) \cdot \nabla v
  = 0 
  \label{primaepsz}
  \\
  \separa
  && \tau \iO \dt\phieps(0) \, v
  + \iO \nabla\phiz \cdot \nabla v
  + \iO \bigl( \betaeps(\phiz) + \pi(\phiz) \bigr) v
  + \iO \seps(\chiz) v
  \qquad
  \non
  \\
  && = \iO \bigl( \mueps(0) + g(0) \bigr) v
  \label{secondaepsz}
\Esist
which hold for every $v\in\Vz$ and every $v\in V$, respectively.
To derive the \pier{needed} estimate for~$\phieps(0)$, we test \eqref{primaepsz} and \eqref{secondaepsz} 
by $\calD\dt\phieps(0)$ and~$\dt\phieps(0)$, respectively.
Then, we get rid of $\mueps(0)$ by adding the equalities we obtain to each other
and owing to the cancellation that occurs due to the definition of~$\calD$.
We have
\Bsist
  && \frac 12 \, \normaVp{\dt\phieps(0)}^2
  + \frac \tau 2 \iO |\dt\phieps(0)|^2
  \non
  \\
  && = - \iO \nabla\phiz \cdot \nabla\dt\phieps(0)
  + \iO \bigl(
    g(0) - \betaeps(\phiz) - \pi(\phiz) - \rho \signeps(\chiz)
  \bigr) \dt\phieps(0) \,.
  \non
\Esist
As for the first term on the \rhs, we recall that $\phiz\in W$ (see~\eqref{hpphizbis}) and~have
\Beq
   - \iO \nabla\phiz \cdot \nabla\dt\phieps(0)
   = \iO \Delta\phiz \, \dt\phieps(0)
   \non
\Eeq
so that we can treat it with the Young inequality.
The same inequality can be used for the second integral we have to estimate
if we recall~that
\Beq
  0 \leq \betaeps(\phiz) \leq \betaz(\phiz) 
  \aand
  |\signeps(\chiz)| \leq 1
  \quad \aeO \,.
  \non
\Eeq 
Hence, the whole \rhs\ is \pier{estimated} from above by
\Bsist
  && \frac \tau 4 \iO |\dt\phieps(0)|^2
  + \frac 1\tau \, \bigl(
    2 \bigl\|
      |\Delta\phiz|
      + |g(0)|
      + |\betaz(\phiz)|
      + |\pi(\phiz)|
    \bigr\|^2
    + 2 \, \normaH{\rho\signeps(\chiz)}^2
  \bigr)
  \non
  \\
  && \leq \frac \tau 4 \iO |\dt\phieps(0)|^2
  + c + \frac {2\mO} \tau \, \rho^2 
  \non
\Esist
and we conclude that
\Beq
  \frac 12 \, \normaVp{\pier{\dt}\phieps(0)}^2
  + \frac \tau 4 \, \normaH{\pier{\dt}\phieps(0)}^2
  \leq  c + \frac {2\mO} \tau \, \rho^2 .
  \non
\Eeq
The next term on the \rhs\ of~\eqref{perquintastima} 
can be treated by means of \eqref{hpgbis} and~\eqref{primastimanuova}, namely
\Beq
  \intQt \dt g \, \dt\phieps
  - \intQt \pi'(\phieps) |\dt\phieps|^2 
  \leq c \, \bigl( 1 + \norma{\dt\phieps}_{\L2H}^2 \bigr)
  \leq c \, (1 + \rho) 
  \non
\Eeq
and the last three terms of \eqref{perquintastima} can be trivially estimated by $c\rho$ 
since $|\seps(r)|\leq\rho$ for every $r\in\erre$.
Hence, by accounting for all the estimate we have obtained,
we deduce from \eqref{perquintastima}~that
\Beq
  \norma{\dt\phieps}_{\L\infty H\cap\L2V}
  \leq C_1 \mO^{1/2} \rho + c (\rho^{1/2} + 1)
  \leq C_2 \mO^{1/2} \rho + c
  \label{quintastima}
\Eeq
where $C_1$ is a constant that depends only on the structure
and, e.g., $C_2:=C_1+1$.

\step
Sixth a priori estimate

We can see \eqref{secondaeps} as a PDE
(see \eqref{pde} in Remark~\ref{Piureg}, which also applies to \eqref{secondaeps}, of course).
We write it at the time~$t$ (\aet),
multiply by $-\Delta\phieps(t)$ and integrate over~$\Omega$.
However, for brevity, we avoid writing the time $t$ for a while.
We obtain
\Bsist
  && \iO |\Delta\phieps|^2
  + \iO \betaeps'(\phieps) |\nabla\phieps|^2
  + \iO \seps'(\chieps) \nabla\chieps \cdot \nabla\phieps
  \non
  \\
  && = \iO \bigl( 
    \mueps + g
    - \tau \dt\phieps
    - \pi(\phieps)    
  \bigr) (-\Delta\phieps) .
  \non
\Esist
We treat the integral involving $\seps'$ \pier{as follows}
\Bsist
  && \iO \seps'(\chieps) \nabla\chieps \cdot \nabla\phieps
  = \iO \seps'(\chieps) |\nabla\chieps|^2
  + \iO \nabla(\seps(\chieps)) \cdot \nabla\phistar
  \non
  \\
  && \geq \iO \nabla(\seps(\chieps)) \cdot \nabla\phistar
  = \iO \seps(\chieps) (-\Delta\phistar)
  \geq - \rho \, \norma{\Delta\phistar}_{\L\infty\Luno}
  = - c \, \rho 
  \non
\Esist
the uniform summability of $\Delta\phistar$ following from \eqref{hpphistarbis}.
We deduce \aat
\Beq
  \frac 12 \, \normaH{\Delta\phieps(t)}^2
  \leq \frac 12 \, \normaH{\mueps(t) + g(t) - \tau \dt\phieps(t) - \pi(\phieps(t))}^2
  + c \, \rho
  \non
\Eeq
whence
\begin{align}
  & \normaH{\Delta\phieps(t)}
  \leq \normaH{\mueps(t) + g(t) - \tau \dt\phieps(t) - \pi(\phieps(t))}
  + c \, \rho^{1/2}
  \non
  \\
  & \leq \normaH{\mueps(t)}
  + \normaH{g(t)}
  + \tau \normaH{\dt\phieps(t)}
  + \normaH{\pi(\phieps(t))}
  + c \, \rho^{1/2} \,.
  \label{pier1}
\end{align}
\pier{We aim to deduce two uniform bounds in~$\LQ\infty$.}
First, from \eqref{primaeps}, by adapting the argument of Remark~\ref{Piureg},
we \pier{have} that $-\Delta\mueps=\dt\phieps$.
Therefore, by \eqref{quintastima} we infer that 
\Beq
  \norma{\Delta\mueps}_{\L\infty H}
  \leq C_2 \mO^{1/2} \rho + c
  \non
\Eeq 
so that the embedding inequality \eqref{embeddingz} yields
\Beq
  \norma\mueps_\infty
  \leq \pier{C_3} \mO^{2/3} \rho + c 
  \label{stimamu}
\Eeq
with $\pier{C_3}:=\Csh C_2$. \pier{Now, by accounting for \eqref{pier1}, \eqref{stimamu}, \eqref{hpgbis},  \eqref{quintastima} and the second estimate in~\eqref{dastimabase}, 
we infer} that 
\begin{align}
  &\norma{\Delta\phieps}_{\pier{\L\infty H}}
  \leq \pier{\mO^{1/2} \left( C_3 \mO^{2/3} +  \tau C_2 \right) \rho + c \, ( \rho^{1/2} + 1 )} \non \\
  &\leq \pier{C_4 \, \mO^{1/2} \left(\mO^{2/3} + 1 \right) \rho + c }
  \label{sestastima}
\end{align}
where, e.g., $\pier{C_4 :=\max\{C_3 , \tau C_2+1\}} $. Next, \pier{in the light of  \eqref{sestastima}, the first condition in \eqref{dastimabase}} and the embedding inequality \eqref{embedding},
we conclude that 
\Beq
  \norma\phieps_\infty
  \leq  \pier{C_5 \, \mO^{2/3} \left(\mO^{2/3} + 1 \right) \rho + c }  \non
\Eeq
where, e.g., $C_5:=\pier{\Csh C_4 } $.
Since $\pi$ is \Lip\ continuous, we deduce that
\Beq
  \norma{\pi(\phieps)}_\infty
  \leq  \pier{C_6 \, \mO^{2/3} \left(\mO^{2/3} + 1 \right) \rho + c }
  \label{stimapi}
\Eeq
with $C_6:=C_5 \sup|\pi'|$. At this point, if we set
\Beq
  \Geps := 
    \mueps + g - \pi(\phieps) - \tau \dt\phistar - \Delta\phistar 
    \label{defGeps}
\Eeq
and recall the estimates \eqref{stimamu} and \eqref{stimapi} as well as our assumptions on $g$ and~$\phistar$, 
we conclude~that
\Beq
  \norma\Geps_\infty
  \leq \Cstr \mO^{2/3}\pier{\left(\mO^{2/3} + 1 \right)} \rho + \hatC
  \label{forcomparison}
\Eeq
where $\Cstr:=\pier{C_3} + C_6$ only depends on the structure of the problem and the shape constant~$\Csh$
(see the construction of the previous~$C_i$'s),
while $\hatC$ also depends on $\Omega$, $T$ and the data $g$, $\phiz$ and~$\phistar$.
\medskip

As already announced, 
we show the existence of a sliding mode by a comparison argument.
The function we use in this project is related to
the solution to a Cauchy problem for an ordinary differential equation
we study at once.

\step
An ordinary differential equation

Given two real numbers $M\in[0,\rho)$ and $\wz\geq0$,
we consider the problem of finding $\weps\in W^{1,\infty}(0,T)$ such~that
\Beq
  \tau \weps'(t) + \rho \signeps \weps(t) = M
  \quad \aat
  \aand
  \weps(0) = \wz \,.
  \label{odeeps}
\Eeq
First of all, since $\signeps$ is \Lip\ continuous, such a problem has a unique solution.
In the (less interesting) case $\wz=0$, the solution is given~by
\Beq
  \weps(t)
  = \frac {\eps M} \rho \Bigl(
    1 - \exp \frac {-\rho\,t}{\eps\tau}
  \Bigr)
  \quad \hbox{for $t\in[0,T]$}.
  \label{wepszero}
\Eeq
Indeed, such a formula provides $\weps$ satisfying 
$\weps(0)=0$ and $\tau\weps'+\rho\weps/\eps=M$.
Since $0\leq\weps\leq\eps M/\rho<\eps$, we also have that $\signeps\weps=\weps/\eps$
so that \eqref{odeeps} is fulfilled.
Suppose now that $\wz>0$.
Then we can assume $\eps\in(0,\wz)$.
We prove that
\Beq
  0 \leq \weps \leq \wz
  \quad \aet \,.
  \label{rangeweps}
\Eeq
In order to derive the first inequality, 
we multiply the equation \eqref{odeeps} by~$-\weps^-$, 
where here and \pier{in the sequel} $(\cdot)^-$ denotes the negative part
(later on, we also use the symbol $(\cdot)^+$ \pier{for} the positive part).
Then, we integrate over~$(0,t)$ and rearrange.
Since $\wz\geq0$, we obtain
\Beq
  \pier{\frac\tau 2} \, |\weps^-(t)|^2
  - \rho \iot (\signeps \weps(s)) \weps^-(s) \ds
  = - M \iot \weps^-(s) \ds
  \leq 0 \,.
\Eeq
On the other hand, $\signeps r\leq 0$ for $r\leq0$
so that the second term on the \lhs\ is nonegative.
Thus $\weps^-=0$, whence $\weps\geq0$.
In order to show the second inequality \pier{in}~\eqref{rangeweps}
we consider the open set $P$ of points $t\in(0,T)$ such that $\weps(t)>\eps$.
We have $\signeps\weps=1$ in~$P$ whence also
$\tau\weps'=M-\rho<0$.
Hence, by recalling that $P$ is~a (finite or countably infinite) union of open intervals~$I_n$,
we infer that the restriction of $\weps$ to each of them is a strictly decreasing affine function.
On the other hand, one of these intervals, say~$I_1$, has $0$~as an end-point since $\weps(0)=\wz>\eps$.
Therefore, one can easily derive that $P=I_1$,
so that there are two possibilities.
It might happen that $\weps\geq\eps$ in the whole of~$[0,T]$.
In this case, $\weps$~is a strictly decreasing affine function.
In the opposite case, $\weps$~is strictly decreasing till it reaches the value~$\eps$ 
at some $T_\eps<T$ and then it remains under such a level.
Thus, the second inequality of \eqref{rangeweps} is established in any case.
We notice that $\weps$~could be explicitly computed,
but the calculation is not necessary.

Finally, we prove that $\weps$ converges as $\eps\searrow0$ to the (unique) solution $w$ to the following problem
\Beq
  \tau w'(t) + \rho \sign w(t) \ni M
  \quad \aat
  \aand
  w(0) = \wz.
  \label{ode}
\Eeq
If $\wz=0$, then $0\leq\weps\leq\eps$, whence $\weps$ tends to zero uniformly.
On the other hand, $w=0$ solves \eqref{ode} since $M\in[0,\rho)\subset\rho\sign 0$.
Suppose now that $\wz>0$. 
Then we can assume $\eps\in(0,\wz)$.
We trivially have $\tau|\weps'|\leq M+\rho$.
Moreover, $\weps(0)$~is independent of~$\eps$.
By also applying the Ascoli--Arzel\`a theorem, we deduce~that
\Beq
  \weps \to w
  \quad \hbox{weakly star in $W^{1,\infty}(0,T)$ and strongly in~$C^0([0,T])$}
  \non
\Eeq
for some function $w\in W^{1,\infty}(0,T)$,
in principle for a subsequence.
Since $\weps(0)$ converges to~$w(0)$, we obtain $w(0)=\wz$.
Next, we recall that $\signeps$ is bounded, 
so that 
\Beq
  \signeps\weps \to \sigma
  \quad \hbox{weakly star in $L^\infty(0,T)$}
  \non
\Eeq
for some $\sigma\in L^\infty(0,T)$
(once more for a subsequence, in principle), whence \pier{we} immediately
\pier{infer that} $\tau w'+\rho\,\sigma=M$ \aet.
By applying, e.g., \cite[Lemma~2.3, p.~38]{Barbu},
we deduce that $\sigma\in\sign w$ so that $w$ solves problem~\eqref{ode}.

\Brem
\label{RemTstar}
The function $w$ can be explicitly computed.
Namely, we have
\Beq
  w(t) = \Bigl( \wz - \frac {\rho-M}\tau \, t \Bigr)^+
  \quad \hbox{for $t\in[0,T]$}.
  \label{soluzode}
\Eeq
In particular, if we define the nonnegative number
\Beq
  \Tstar := \frac {\tau \wz}{\rho-M}
  \label{defTstar}
\Eeq
and we reinforce our assumption on $\rho$ by assuming that 
\Beq
  \rho > M + \frac {\tau \wz} T
  \label{condrho}
\Eeq
then $\Tstar<T$ and $w(t)=0$ for every $t\in[\Tstar,T]$.
\Erem

\step
The comparison argument

The function we use in our argument is the space independent function
$(x,t)\mapsto\weps(t)$ (still termed $\weps$ for simplicity)
where $\weps$ is the solution to \eqref{odeeps} with a proper choice of $\wz$ and~$M$.
We recall that $\Cstr$ and $\hatC$ are the \pier{constants} that appear in~\eqref{forcomparison}.
We stress once more that $\Cstr$ only depends on the structure of the problem and the shape constant~$\Csh$.
We assume that
\Beq
  \mO < \deltastar  := \pier{\left|\frac{ \sqrt{1+ 4/\Cstr } - 1}2 \right|^{3/2}}  
  \label{hpmO}
\Eeq
\pier{so that} $ \Cstr \mO^{2/3}\pier{\left(\mO^{2/3} + 1 \right)} <1$.
By also recalling \eqref{hpphistarbis},
we see that the real number
\Beq
  \rhostar := \frac {\hatC + \betastar + (\tau\wz)/T} {1 - \Cstr \mO^{2/3}\pier{\left(\mO^{2/3} + 1 \right)}}
  \quad \hbox{where} \quad
  \betastar := \norma{\betaz(\phistar)}_\infty
  \label{defrhostar}
\Eeq 
is \pier{well defined}. At this point, we fix $\rho$ by assuming~that 
\Beq
  \rho > \rhostar 
  \label{hprho}
\Eeq
and choose
\Beq
  \wz := \norma\chiz_\infty 
  \aand
  M := \Cstr \mO^{2/3} \pier{\left(\mO^{2/3} + 1 \right)}\rho + \hatC + \betastar \,.
  \label{scelte}
\Eeq
Our assumptions and choices are made in order that $\rho>M+(\tau\wz)/T$, whence in particular $\rho>M$.
Hence, the conditions assumed in the study of the solution $\weps$ performed above 
and in Remark~\ref{RemTstar} are fulfilled,
and the time $\Tstar$ given by \eqref{defTstar} is well-defined and belongs to~$(0,T)$.
At this point, we can start our comparison argument.
We recall that $\phistar(t)\in W$ \aat\ by \eqref{hpphistarbis} 
so that we can integrate by parts in the \rhs\ of \eqref{secondaepsbis}
and write the equation~as
\Bsist
  && \tau \iO \dt\chieps \, v
  + \iO \nabla\chieps \cdot \nabla v
  + \iO \betaeps(\chieps+\phistar)
  + \rho \iO \signeps(\chieps) v
  = \iO \Geps v
  \qquad
  \non
  \\
  && \quad \hbox{\aet\ and for every $v\in V$}
  \label{secondaepster}
\Esist
where $\Geps$ is given by~\eqref{defGeps}.
We recall \eqref{forcomparison}, which provides a uniform bound for~$\Geps$.
On the other hand, by reading $\weps$ as a space independent function defined in~$Q$ as said before,
we can write the ordinary differential equation \eqref{odeeps}
as a partial differential equation with homogeneous Neumann boundary conditions.
It is convenient to choose the following two forms
\Bsist
  && \tau \iO \dt\weps \, v
  + \iO \nabla\weps \cdot \nabla v
  + \iO \bigl(
    \betaeps(\pier{\weps} +\phistar) + \rho \signeps\weps
  \bigr) v
  \non
  \\
  && = \iO \bigl(
    M + \betaeps(\pier{\weps}+\phistar)
  \bigr) v
  \quad \hbox{\aet\ and for every $v\in V$}
  \label{persopra}
  \\
  && \tau \iO \dt\weps \, v
  + \iO \nabla\weps \cdot \nabla v
  + \iO \bigl(
    -\betaeps(-{}\pier{\weps{}}+\phistar) - \rho \signeps(-\weps)
  \bigr) v
  \non
  \\
  && 
  = \iO \bigl(
    M - \betaeps(-\pier{{}\weps{}}+\phistar)
  \bigr) v
  \quad \hbox{\aet\ and for every $v\in V$} .
  \label{persotto}
\Esist
We prove that $|\chieps|\leq\weps$ \aeQ\ by showing that
\Beq
  \chieps \leq \weps
  \aand
  - \chieps \leq \weps
  \quad \aeQ \,.
  \label{tesi}
\Eeq
To obtain the first inequality, we take the difference between \eqref{secondaepster} and \eqref{persopra}
and \pier{choose $v= (\chieps-\weps)^+$}.
Then, we integrate over~$(0,t)$.
We~have
\Bsist
  && \frac \tau 2 \iO |(\chieps-\weps)^+(t)|^2
  + \intQt |\nabla(\chieps-\weps)^+|^2
  \non
  \\
  && \quad {}
  + \intQt \{
    \betaeps(\chieps+\phistar) - \betaeps(\weps+\phistar)
    + \rho \signeps\chieps - \rho \signeps\weps
  \} (\chieps-\weps)^+
  \non
  \\
  && = \intQt \{
    \Geps - M - \betaeps(\weps+\phistar)
  \} (\chieps-\weps)^+ \,.
  \non
\Esist
Clearly, the expression between braces on the \lhs\ is nonnegative in the set where $\chieps>\weps$
so that the corresponding integral is nonnegative.
On the other hand, since $\weps$ is nonnegative and recalling the estimate \eqref{forcomparison}
and the \pier{definitions} of $\betastar$ and~$M$ (see \eqref{defrhostar} and~\eqref{scelte}), we~have
\Beq
  \Geps - M - \betaeps(\weps+\phistar)
  \leq \norma\Geps_\infty - M - \betaeps(\phistar)
  \leq \norma\Geps_\infty - M + |\betaz(\phistar)|
  \leq 0 \,.
  \non
\Eeq
Hence, $(\chieps-\weps)^+=0$ and the desired inequality is proved.
To obtain the other one, we add equations \eqref{secondaepster} and \eqref{persotto} to each other
and test the equality we get by~$-(\chieps+\weps)^-$.
Then, we integrate over~$(0,t)$.
We~have
\Bsist
  && \frac \tau 2 \iO |(\chieps+\weps)^-(t)|^2
  + \intQt |\nabla(\chieps+\weps)^-|^2
  \non
  \\
  && \quad {}
  - \intQt \{
    \betaeps(\chieps+\phistar) - \betaeps(-\weps+\phistar)
    + \rho \signeps\chieps - \rho \signeps(-\weps)
  \} (\chieps+\weps)^-
  \non
  \\
  && = - \intQt \{
    \Geps + M - \betaeps(-\weps+\phistar)
  \} (\chieps+\weps)^- \,.
  \non
\Esist
In the set where $\chieps+\weps$ is negative, we have $\chieps<-\weps$
so that the expression between braces on the \lhs\ is nonpositive 
and the corresponding integral is nonpositive.
On the other hand, we have
\Beq
  \Geps + M - \betaeps(-\weps+\phistar)
  \geq - \norma\Geps_\infty + M - \betaeps(\phistar)
  \geq - \norma\Geps_\infty + M - |\betaz(\phistar)|
  \geq 0 \,.
  \non
\Eeq
Hence $(\chieps+\weps)^-=0$ and \eqref{tesi} is completely proved.
Since $\chieps$ and $\weps$ converge to $\chi$ and~$w$, respectively,
where \pier{$w$ is the solution to~\eqref{ode}},
we deduce that
\Beq
  |\chi| \leq w 
  \quad \aeQ \,.
  \non
\Eeq
As already noticed,
we can apply Remark~\ref{RemTstar}.
Hence, $\Tstar<T$ and $\chi(t)=0$ for every $t\in[\Tstar,T]$, i.e.,
$\phi(t)=\phistar(t)$ for every $t\in[\Tstar,T]$.

\gabri{%
\Brem
In terms of the original physical variables,
i.e., the order parameter $\phi$ and the associated chemical potential $\mu$ satisfying the inhomogeneous Dirichlet boundary condition,
the behavior of the solution after the time $\Tstar$ is the following:
\Bsist
  && \phi(t) = \phistar(t)
  \quad \hbox{for every $t\in[\Tstar,T]$}
  \non
  \\
  && -\Delta\mu(t) = -\dt\phistar(t)
  \aand
  \mu(t)|_\Gamma = \muG(t)
  \quad \hbox{for a.a.\ $t\in(\Tstar,T)$}.
  \non
\Esist
In particular, even $\mu$ can be explicitly computed on~$(\Tstar,T)$.
\Erem
}%

\Brem
Due to our choice \eqref{scelte} of~$M$, the difference $\rho-M$ is almost proportional to~$\rho$ for large values of~it.
It follows that the time $\Tstar$ given by~\eqref{defTstar} tends to zero as $\rho$ tends to infinity.
Therefore, the sliding mode can be imposed to occur in an arbitrarily short time 
by assuming that $\rho$ is large enough.
\Erem


\section*{Acknowledgements} 

The authors warmly thank Professor
Viorel Barbu for some useful discussions and suggestions. This research
activity has been performed in the framework of the Italian-Romanian
collaboration agreement \textquotedblleft\pier{Analysis and optimization of 
mathematical models ranging from bio-medicine to
engineering}\textquotedblright\ between the
Italian CNR and the Romanian Academy. \oldpier{PC acknowledges support 
from the Italian Ministry of Education, 
University and Research~(MIUR): Dipartimenti di Eccellenza Program (2018--2022) 
-- Dept.~of Mathematics ``F.~Casorati'', University of Pavia, and from 
the GNAMPA (Gruppo Nazionale per l'Analisi Matematica, 
la Probabilit\`a e le loro Applicazioni) of INdAM (Isti\-tuto 
Nazionale di Alta Matematica).} The present paper also benefits from the
support \oldpier{of the Romanian Ministry of Research} and Innovation, CNCS --UEFISCDI,
project number PN-III-P4-ID-PCE-2016-0011, for~GM.



\end{document}